\documentstyle[11pt]{article}

\begin{document}
\title{\Large Linear maps between $C^\ast$-algebras whose adjoints preserve 
		extreme points of the dual ball\vspace{1cm}}
  
\author{
   Louis E Labuschagne\\
   {\small Department of Mathematics \& Applied Mathematics,
   University of Pretoria}\\
   {\small 0002 Pretoria, South Africa $\quad$
   [llabusch@scientia.up.ac.za]}
\and
    Vania Mascioni\\
    {\small Department of Mathematics, The University of Texas at Austin}\\
    {\small Austin, TX 78712, U.S.A. $\quad$
    [mascioni@math.utexas.edu]}
    }
\date{February 1996}
\maketitle

\newpage
\vspace*{5cm}

Short running title: Extremal maps on $C^\ast$-algebras

Mailing Address : 

    \begin{verse}
    Vania Mascioni\\
    Department of Mathematics\\
    The University of Texas at Austin\\
    Austin, TX 78712, U.S.A.\\    
    e-mail: mascioni@math.utexas.edu\\
    \end{verse}

\newpage
\vspace*{5cm}

\begin{abstract}
We give a structural characterisation of linear operators from one $C^\ast$%
-algebra into another whose adjoints map extreme points of the dual ball
onto extreme points. We show that up to a $\ast$-isomorphism, such a map
admits of a decomposition into a degenerate and a non-degenerate part, the
non-degenerate part of which appears as a Jordan $\ast$-morphism followed by
a ``rotation'' and then a reduction. In the case of maps whose adjoints
preserve pure states, the degenerate part does not appear, and the
``rotation'' is but the identity. In this context the results concerning such
pure state preserving maps depend on and complement those of St\o rmer [St\o
2; 5.6 \& 5.7]. In conclusion we consider the action of maps with ``extreme
point preserving'' adjoints on some specific $C^\ast$-algebras.
\end{abstract}
\newpage

\section{Introduction}

It is clear from the remarks made in the abstract that the results
concerning maps with ``pure state preserving'' adjoints, provide us with a
valuable clue as to what objects we may regard as ``non-commutative
composition operators''. The value of these and the other results also lie in
the fact that they indicate that results of this nature for $C(K)$ spaces
are not merely isolated fragments, but rather indicative of a very deep $%
C^\ast$-algebraic structure reaching far beyond the simplicity of the
commutative case.

The notation employed is fairly standard $C^\ast$-algebraic notation and for
the most part amounts to a subtle interpolation of that of Bratteli and
Robinson [BR], and Kadison and Ringrose [KR]. The main features are the
following:

${\cal A, B}$ and ${\cal C}$ will be deemed to be typical $C^\ast$%
-algebras which for the sake of convenience we will assume to be unital.
Given ${\cal A}$, the associated sets of all states and all pure states
of ${\cal A}$ will be denoted by ${\cal S_A}$ and ${\cal P_A}$
respectively. If indeed ${\cal A}$ is concrete, ${\cal A^{\prime}}$
denotes the commutant and ${\cal N_A}$ the set of all normal states.
Functionals of a $C^\ast$-algebra will be denoted by $\rho, \omega$, with $%
\omega$ being reserved for the notation of states (usually pure). In this
context, given a state $\omega$ of ${\cal A},$ \ $(\pi_\omega, h_\omega,
\Omega_\omega)$ will denote the canonical cyclic representation of ${\cal %
A}$ engendered by $\omega$. Here $h_\omega$ is the relevant Hilbert space, $%
\Omega_\omega \in h_\omega$ the cyclic unit vector corresponding to $\omega$%
, and $\pi_\omega$ the canonical $\ast$-homomorphism from ${\cal A}$ into 
$B(h_\omega)$. Typical Hilbert spaces will be taken to be $h$ and $k$.
Finally given any Banach space $X, \ (X)_1$ or $X_1$ if there is no danger
of confusion, will denote the closed unit ball of $X$. In this context $%
{\rm ext}(X_1)$ denotes the set of extreme points of $X_1$.

Regarding linear maps on $C^{*}$-algebras, a Jordan ($*$-)morphism is
understood to be a mapping $\psi :{\cal A}\rightarrow {\cal B}$ such
that $\psi (AB+BA)=\psi (A)\psi (B)+\psi (B)\psi (A)$ and $\psi (A^{*})=\psi
(A)^{*}$ for all $A,B\in {\cal A}$. This concept is of course equivalent
to that of a $C^{*}$-homomorphism which is defined to be a positive map
preserving squares of self-adjoint elements. To see this one need only note
that in general $(A+B)^2-A^2-B^2=AB+BA$, and make use of the fact ${\rm %
span}({\cal A}_{sa})={\cal A}$. Moreover given any Jordan $*$-morphism 
$\psi :{\cal A}\rightarrow {\cal B},\ \psi (I)=E$ is easily seen to be
an orthogonal projection with $\psi (A)=E\psi (A)E$ for all $A\in {\cal A}
$. To see the latter fact one need only note that if indeed $\psi $ is a
Jordan $*$-morphism, then $\psi (ABA)=\psi (A)\psi (B)\psi (A)$ for all $%
A,B\in {\cal A}$ [BR; p 212]. Thus as a map into ${\cal B}_E$, \ $\psi 
$ then preserves the identity. In particular if ${\cal B}$ is concrete
and $\psi $ a Jordan $*$-morphism with $\psi ({\cal A})^{\prime \prime }=%
{\cal B}^{\prime \prime }$, we must then have $\psi (I)=I$. In this
context we also observe that for our purposes we do not need to assume
continuity of the operators we characterise, since the properties under
consideration necessarily imply that these must even have norm one. In the
case where $\psi :{\cal A}\rightarrow {\cal B}$ with $\omega \circ
\psi \in {\cal P_A}$ for every $\omega \in {\cal P_B}$, this follows
from [BR; 3.2.6] on noticing that by [KR; 4.3.8] we have $\psi \geq 0$ with $%
\psi (I)=I$. In the case where $\rho \circ \psi \in {\rm ext}({\cal A}%
_1^{*})$ for every $\rho \in {\rm ext}({\cal B}_1^{*})$, we merely
need to verify continuity and apply the Krein-Milman theorem to $\psi ^{*}$.
To see continuity in this case, given $A\in {\cal A}$, select $\omega
_0,\omega _1\in {\cal P_B}$ so that $\omega _0(\mbox{\rm Re}\,\psi
(A)) = \Vert \mbox{\rm Re}\,\psi (A)\Vert $ and $\omega _1(\mbox{\rm Im}\,
\psi (A))=\Vert \mbox{\rm Im}\,\psi (A)\Vert $ [KR; 4.3.8]. Then since $\omega
_0\circ \psi ,\omega _1\circ \psi \in {\rm ext}({\cal A}_1^{*})$, they
are both norm-one functionals and hence

\begin{eqnarray*}
\Vert \psi (A)\Vert & \leq & \Vert \mbox{\rm Re}\,\psi (A)\Vert +
  \Vert \mbox{\rm Im}\,\psi (A)\Vert \\ 
& = & \omega _0(\frac 12(\psi (A)+\psi (A)^{*}))+\omega _1(\frac{-i}2(\psi
(A)-\psi (A)^{*})) \\ 
& = & \frac 12[\omega _0(\psi (A))+\overline{\omega _0(\psi (A))}-i\omega
_1(\psi (A))+i\overline{\omega _1(\psi (A))}] \\ 
& \leq & (\Vert \omega _0\circ \psi \Vert +\Vert \omega _1\circ \psi \Vert
)\Vert A\Vert \\ 
& = & 2\Vert A\Vert
\end{eqnarray*}
as required.

\section{Maps with pure state preserving adjoints: The overture to the
general case}

Although the lemmas in this section may be deemed to be standard folklore
and Theorem 5 judged to be a technical reworking of hard work done by St\o
rmer ([St\o 1], [St\o 2]), its value lies in the fact that it does present a
coherent framework within which to attack the more general case of maps
whose adjoints preserve extreme points of the unit ball. (To see that this
case is indeed more general is none too trivial (cf. Corollary 20).) Like
any good overture, this section and its lemmas presents in embryonic form
the main ideas developed further later on. For this reason we have chosen to
prove the lemmas in full.

{\bf Lemma 1} \quad Let ${\cal A}$ be a $C^{*}$-algebra and $E\in 
{\cal A}$ a projection. Denote the reduction ${\cal A}\rightarrow 
{\cal A}_E$ by $\eta $. Then $\omega \circ \eta $ is a pure state of $%
{\cal A}$ whenever $\omega $ is a pure state of ${\cal A}_E.$
Conversely if $\tilde{\omega}(E)=1$, then the restriction of $\tilde{\omega}$
to ${\cal A}_E$ is a pure state of ${\cal A}_E$ whenever $\tilde{\omega%
}$ is a pure state of ${\cal A}$.

{\bf Proof} \quad Suppose $\omega $ is a pure state of ${\cal A}_E$
and let $\rho $ be a positive functional on ${\cal A}$ majorised by $%
(\omega \circ \eta )$. We show that then $\rho (A)=\rho (EAE)$ for every $%
A\in {\cal A}$. If this be true, then clearly $\rho $ is of the form $%
\rho _E\circ \eta $ where $\rho _E$ is the restriction of $\rho $ to $%
{\cal A}_E$. Since $\omega \circ \eta \geq \rho =\rho _E\circ \eta \geq 0$%
, it is clear that then $\omega \geq \rho _E\geq 0$. But then $\rho _E$ will
be a multiple of $\omega $ on ${\cal A}_E$ [KR; 3.4.6] and hence $\rho
=\rho _E\circ \eta $ a multiple of $\omega \circ \eta $. By [KR; 3.4.6], $%
\omega \circ \eta $ must then be pure. In order to finally verify that $\rho
(A)=\rho (EAE)$ for every $A\in {\cal A}$, it suffices to do this for the
case $A\in {\cal A}^{+}$ since ${\rm span}({\cal A}^{+})={\cal A}
$. Now if $A\in {\cal A}^{+}$ and $0\leq \rho \leq \omega \circ \eta $,
then surely $0\leq (I-E)A(I-E)$, and hence 

\[
0\leq \rho ((I-E)A(I-E))\leq (\omega \circ \eta )((I-E))A(I-E))=\omega (0)=0
\]

for every $A\in {\cal A}^{+}$, that is

\begin{equation}
\rho((I-E)A(I-E)) = 0.
\end{equation}

Next appealing to (1) and applying [KR; 4.3.1], we get

\begin{eqnarray}
|\rho((I - E)AE)|^2 & \leq & \rho(E^\ast E) \rho((I - E)A((I - E)A)^\ast) \\
& = & \rho(E) \rho((I - E)A^2(I - E)) = 0  \nonumber
\end{eqnarray}

for every $A \in {\cal A}^+$, and hence also that

\begin{eqnarray}
|\rho(EA(I - E))| & = & |\overline{\rho(EA(I - E))}| = |\rho((EA(I -
E))^\ast )| \\
& = & |\rho((I - E)AE)| = 0.  \nonumber
\end{eqnarray}

Combining (1), (2) and (3), we have $\rho(A) = \rho(EAE)$ for every $A \in 
{\cal A}^+$ as required.

Conversely if $\tilde{\omega} \in {\cal P}_{{\cal A}}$ with $\tilde{%
\omega}(E) = 1$, then it may be verified that $\tilde{\omega}(A) = \tilde{%
\omega}(EAE)$ for any $A \in {\cal A}$. As before it suffices to show
this for the case $A \in {\cal A}^+$. For $A \in {\cal A}^+$ it may
easily be verified that

\begin{equation}
0 = \tilde{\omega}((I - E)AE) = \tilde{\omega}(EA(I - E)) = \tilde{\omega}%
((I - E)A(I - E)).
\end{equation}

We show how to do this in one of the cases, the others being similar. Since $%
\tilde{\omega}(I - E) = \tilde{\omega}(I) - \tilde{\omega}(E) = 0$, we have
by [KR; 4.3.1] that

\[
\begin{array}{lll}
|\tilde{\omega}((I - E)AE)|^2 & \leq & \tilde{\omega}((I - E)(I - E)^\ast) 
\tilde{\omega}((AE)^\ast(AE)) \\ 
& = & \tilde{\omega}(I - E) \tilde{\omega}(EA^2E) = 0
\end{array}
\]

for every $A \in {\cal A}^+$. Thus $\tilde{\omega}(A) = \tilde{\omega}%
(EAE)$ for all $A \in {\cal A}^+$ by (4), as required. Clearly then $%
\tilde{\omega}$ is of the form $\omega_0 \circ \eta$ where $\omega_0$ is the
restriction of $\tilde{\omega}$ to ${\cal A}_E.$ Moreover $\omega_0$ is a
state of ${\cal A}_E$ by [KR; 4.3.2] applied to the fact that $\tilde{%
\omega}(E) = 1$. Now finally if $\omega_0 \geq \rho \geq 0$ for some
functional on ${\cal A}_E$, then since $\eta$ preserves order [KR;
4.2.7], we have $\tilde{\omega} = \omega_0 \circ \eta \geq \rho \circ \eta
\geq 0$. Since $\tilde{\omega} \in {\cal P_A}, \ \rho \circ \eta$ must be
a multiple of $\tilde{\omega}$ [KR; 3.4.6] and hence on restriction to $%
{\cal A}_E, \ \rho$ must then be a multiple of $\omega_0$. It follows
that $\omega_0$ is a pure state of ${\cal A}_E$ [KR; 3.4.6]. \hfill $\Box$

{\bf Lemma 2} \quad Let ${\cal A}$ be a von Neumann algebra, $E$ a
projection in ${\cal A}$, and let $\eta$ be defined as before. Then $\rho
\circ \eta$ is a normal state on ${\cal A}$ whenever $\rho$ is a normal
state on ${\cal A}_E$. Conversely if $\tilde{\rho}$ is a normal state of $%
{\cal A}$ with $\tilde{\rho}(E) = 1$, then the restriction of $\rho$ to $%
{\cal A}_E$ is a normal state of ${\cal A}_E$.

{\bf Proof} \quad For the second part all we need to do is note that the
restriction is a state by [KR; 4.2.3], and apply the definition [KR;
7.1.11]. To see the first part all we really need to do is to note that if $%
A_\lambda $ is a monotone increasing net in ${\cal A}$ with least upper
bound $A\in {\cal A}$, then $EA_\lambda E$ is a monotone increasing net
with l.u.b. $EAE$. Then, we use the fact that $\eta $ preserves order, and a
combination of [BR; 2.4.1 \& 2.4.19]. \hfill $\Box $

{\bf Lemma 3} \quad Let ${\cal A}$ be a $C^\ast$-algebra. If ${\cal %
A}$ is in its reduced atomic representation, then every pure state of $%
{\cal A}$ is normal (ultra-weakly continuous).

{\bf Proof} \quad By [KR; 7.1.12], it suffices to show that all the pure
states are vector states. First of all by the definition of the reduced
atomic representation there is a maximal disjoint set of pure states, $%
{\cal M}$, in terms of which the representation is generated by the GNS
construction. These pure states are then obviously vector states. Next given
any $\omega \in {\cal P_A}$, by the maximality of ${\cal M}$, \ $%
\omega $ is unitarily equivalent to some $\omega _0\in {\cal M}$ [KR;
10.2.6 \& 10.3.7], say 
\[
\omega (A)=\omega _0(U^{*}AU)\quad A\in {\cal A}
\]
where $U\in {\cal A}$ is unitary. But then if $\omega _0\in {\cal M}$
corresponds to the vector state say $(A\Omega _0,\Omega _0),\ A\in {\cal A%
},\ \Vert \Omega _0\Vert =1,$ then surely $\omega $ corresponds to 
\[
\omega (A)=\left\langle U^{*}AU\Omega _0,\Omega _0\right\rangle
=\left\langle A(U\Omega _0),(U\Omega _0)\right\rangle \quad A\in {\cal A}
\]
where $\Vert U\Omega _0\Vert =\Vert \Omega _0\Vert =1$. \hfill $\Box $

{\bf Lemma 4} \quad If $\omega$ is a pure normal state of a concrete $%
C^\ast$-algebra ${\cal A}$, the unique normal extension of $\omega$ to $%
{\cal A}^{\prime\prime}$, say $\tilde{\omega}$, is a pure state of $%
{\cal A}^{\prime\prime}$.

{\bf Proof} \quad It is an exercise to show that the ultra-weak
continuity of $\omega $ implies that the representation of ${\cal A}$
engendered by $\omega $ is similarly continuous. (This can be seen by for
example suitably adapting the first part of the proof of [BR; 2.4.24].) If $%
\pi _\omega $ is this representation, then by [KR; 10.1.10] it has a unique
ultra-weakly continuous extension $\tilde{\pi}_\omega $ to all of ${\cal A%
}^{\prime \prime }$. If now $\omega $ corresponds to the vector state $%
\omega _\Omega $ in the sense that $\omega =\omega _\Omega \circ \pi _\omega 
$, then surely $\omega _\Omega \circ \tilde{\pi}_\omega $ is a normal
(ultra-weak) extension of $\omega$, and hence by the uniqueness of
this extension we have $\tilde{\omega}=\omega _\Omega \circ \tilde{\pi}%
_\omega $. A combination of [KR; 10.2.3] and [KR; 10.2.5] applied to $\pi
_\omega $ and $\tilde{\pi}_\omega $ respectively, reveal that 
$\tilde{\omega}$ is a
pure state of ${\cal A}^{\prime\prime}$. \hfill $\Box $

{\bf Theorem 5} \quad Let ${\cal A}$ and ${\cal B}$ be $C^{*}$%
-algebras and $\alpha $ the reduced atomic representation of ${\cal B}$.
Then a linear mapping $\varphi :{\cal A}\rightarrow {\cal B}$ has the
property that $\omega \circ \varphi $ is a pure state whenever $\omega $ is
a pure state of ${\cal B}$ if and only if there exists a von Neumann
algebra ${\cal R}$ acting on some Hilbert space $h$, a projection $E\in 
{\cal R}$, and a set $(F_\nu )$ of mutually orthogonal central
projections in ${\cal R}$ with $\sum_\nu F_\nu =I_{{\cal R}}$, such
that up to a (ultra-weakly continuous) $*$-isomorphic embedding $\Phi $ of $%
\alpha ({\cal B})^{\prime \prime }$ in ${\cal R}$, $\alpha ({\cal B})^{\prime
\prime }$ appears as ${\cal R}_E$ with $\Phi \circ \alpha \circ \varphi $
of the form 
\[
(\Phi \circ \alpha \circ \varphi )(A)=E\psi (A)E\quad\mbox{\rm for all}
    \quad A\in {\cal A}.
\]

Here $\psi $ is a Jordan $*$-morphism from ${\cal A}$ into ${\cal R}$
with the pro\-per\-ty that 
\[
(F_\nu \psi ({\cal A})F_\nu )^{\prime \prime }=F_\nu {\cal R}F_\nu 
\]
(a slightly weaker cond\-i\-tion than me\-rely re\-qui\-ring $\psi ({\cal %
A})^{\prime \prime }={\cal R}$).

{\bf Proof} \quad First assume $\varphi $ to be of the form described in
the hypothesis. Since $*$-isomorphisms clearly preserve pure states, $\Phi
\circ \alpha $ basically identifies ${\cal B}$ with $\Phi (\alpha (%
{\cal B}))$ as far as we are concerned, and hence we may regard ${\cal %
B}$ as a subalgebra of ${\cal R}$ with the property that ${\cal B}%
^{\prime \prime }={\cal R}_E$. We proceed to show that $\varphi $
preserves pure states. Let $\omega \in {\cal P_B}$ be given. By Lemmas 3
and 4 there exists a unique extension $\tilde{\omega}$ of $\omega $ to all
of ${\cal B}^{\prime \prime }={\cal R}_E$ which is pure and
ultra-weakly continuous on ${\cal B}^{\prime \prime }$. By the uniqueness
we may identify $\omega $ with $\tilde{\omega}$. Considering Lemma 1, it
follows that $\omega _E=\omega (E\cdot E)$ is a pure state on ${\cal R}$.
Therefore by [KR; 4.3.14] $\omega _E(F_\nu )\in \{0,1\}$ for every $\nu $.
However since the $F_\nu $'s are mutually orthogonal central projections
with $\vee _\nu F_\nu =I$ and since $\omega _E$ is suitably continuous by
Lemma 3, it follows that $1=\omega _E(I)=\sum_\nu \omega _E(F_\nu )$ and
hence that $\omega _E(F_\nu )=1$ for precisely one of the $F_\nu $'s, say $%
\omega _E(F_{\nu _0})=1.$ Thus denoting $\omega _E(F_{\nu _0}\cdot F_{\nu
_0})$ and $F_{\nu _0}\psi F_{\nu _0}$ by $\omega _{\nu _0}$ and $\psi _{\nu
_0}$ respectively, it follows from [KR; 4.3.14] that 
\[
\omega _E\circ \psi =\omega _{\nu _0}\circ \psi _{\nu _0}.
\]
Clearly it suffices to consider $\omega _E$ in terms of the von Neumann
algebra ${\cal R}_{F_{\nu _0}}=(F_{\nu _0}\psi ({\cal A})F_{\nu
_0})^{\prime \prime }$ only. By Lemmas 1 and 2, $\omega _{\nu _0}$ does
indeed define an ultra-weakly continuous pure state on ${\cal R}_{F_{\nu
_0}}$. Thus we have reduced matters to the case where we have a von Neumann
algebra ${\cal R}_0$, a Jordan-morphism $\psi _0:{\cal A}\rightarrow 
{\cal R}_0$ with the property that the $C^{*}$-algebra ${\cal C}$
generated by $\psi _0({\cal A})$ has ${\cal R}_0$ as its double
commutant, and an ultra-weakly continuous pure state $\omega _0$ on $%
{\cal R}_0$. Now assume that ${\cal C}$ is its own universal
representation, and hence that ${\cal R}_0$ is the bidual of ${\cal C}$. 
If this was not the case we could have ``lifted'' the original description
to this case by means of an application of [KR; 10.1.12] combined with
Lemmas 1 and 2. Since $\omega _0$ is both pure and normal, it is an extreme
point of the set of normal states. Finally by combining for example [KR;
7.4.2, 10.1.1 \& 10.1.2], the set of normal states on ${\cal R}_0$ is
isometrically isomorphic to the state space of ${\cal C}$ under
restriction to ${\cal C}$. Hence the restriction of $\omega _0$ to $%
{\cal C}$ is a pure state of ${\cal C}$. On applying [St\o 1,
Corollary 5.8], we conclude that $(\omega _0\circ \psi _0)$ is a pure state
of ${\cal A}$ as required.

For the converse we first show that for any pure state $\omega$ acting on $%
{\cal B}, \ \pi_\omega \circ \varphi$ has the required form where $%
\pi_\omega$ corresponds to the canonical irreducible representation
generated by $\omega$ [KR; 10.2.3], before deducing the result from this
fact. This is basically a straightforward consequence of [St\o 2, Thm 5.7].
Given $\omega \in {\cal P_B}$, we consider two cases:

If $\pi _\omega \circ \varphi $ is a pure state on ${\cal A}$, then on
denoting $\pi _\omega \circ \varphi $ by $\rho $, let $\pi _\rho $ be the
irreducible representation of ${\cal A}$ on a Hilbert space $h_\rho $
with cyclic unit vector $\Omega _\rho $, generated by $\rho $ by means of
the GNS process. Since $\pi _\rho ({\cal A})$ is irreducible, $\pi _\rho (%
{\cal A})^{\prime \prime }=B(h_\rho )$, and hence the orthogonal
projection $E_\rho $ of $h_\rho $ onto the ray ${\rm span}\{\Omega _\rho
\}$, belongs to $\pi _\rho ({\cal A})^{\prime \prime }$. Since now $%
E_\rho $ is of the form $E_\rho a=\left\langle a,\Omega _\rho \right\rangle
\Omega _\rho $ for any $a\in h_\rho $, it follows that 
\[
\begin{array}{lll}
\left\langle E_\rho \pi _\rho (A)E_\rho a,b\right\rangle & = & \left\langle
\pi _\rho (A)E_\rho a,E_\rho b\right\rangle \\ 
& = & \left\langle \pi _\rho (A)\left\langle a,\Omega _\rho \right\rangle
\Omega _\rho ,\left\langle b,\Omega _\rho \right\rangle \Omega _\rho
\right\rangle \\ 
& = & \left\langle a,\Omega _\rho \right\rangle \overline{\left\langle 
b,\Omega _\rho \right\rangle} \left\langle \pi _\rho (A)\Omega _\rho ,\Omega
_\rho \right\rangle \\ 
& = & \left\langle a,\Omega _\rho \right\rangle \overline{\left\langle
b,\Omega _\rho \right\rangle }\rho (A)\cdot \Vert \Omega _\rho \Vert ^2 \\ 
& = & \rho (A)\left\langle \left\langle a,\Omega _\rho \right\rangle \Omega
_\rho ,\left\langle b,\Omega _\rho \right\rangle \Omega _\rho \right\rangle
\\ 
& = & \rho (A)\left\langle E_\rho a,E_\rho b\right\rangle \\ 
& = & \left\langle \rho (A)E_\rho a,b\right\rangle
\end{array}
\]

for all $A\in {\cal A}$ and all $a,b\in h_\rho $.

Hence $\rho (A)E_\rho =E_\rho \pi _\rho (A)E_\rho $ for every $A\in {\cal %
A}$. In the obvious way we may now identify $\rho ({\cal A})$ with $\rho (%
{\cal A})E_\rho =E_\rho \pi _\rho ({\cal A})E_\rho $, from which it
now follows that $\pi _\omega \circ \varphi =\rho $ is of the required form.

Next suppose $\omega \in {\cal P_B}$, but that $\pi_\omega \circ \varphi$
is not a pure state. In the notation of [St\o 2, Thm 5.7] it now follows that

\[
\pi_\omega \circ \varphi = V^\ast \rho V
\]

where $\rho $ is a Jordan $*$-morphism with $\rho ({\cal A})^{\prime
\prime }=B(h)$ (i.e. the $C^{*}$-algebra generated by $\rho ({\cal A})$
is irreducible), and $V$ is a linear isometry from $h_\omega $ into $h$. Let 
$E\in {\cal B}(h)$ be the orthogonal projection onto $V(h_\omega )$.
Since $VV^{*}=E$ and since $V^{*}\mid _{E(h)}$ is the partial inverse of $V$
on $E(h)$ with $V^{*}=V^{*}\mid _{E(h)}\cdot E$ and $EV=V$, it follows that $%
V$ generates a spatial $*$-isomorphism $\Phi _\omega $ from $B(Eh)$ onto $%
B(h_\omega )$ such that $\Phi _\omega (EB(h)E)=V^{*}B(h)V=B(h_\omega )=\pi
_\omega ({\cal B})^{\prime \prime }$. Clearly we may therefore assume $%
\pi _\omega \circ \varphi $ to be of the form $E\rho E$ as required.

On applying Zorn's lemma, we may now select a maximal set of pure states $%
{\cal M}\subset {\cal P_B}$ such that the irreducible [KR; 10.2.3] GNS
representations generated by any two elements of ${\cal M}$ are mutually
disjoint (pairwise inequivalent) [KR; 10.3.7]. Then surely $\pi
=\bigoplus_{\omega \in {\cal M}}\pi _\omega $ is faithful with 
\[
\pi ({\cal B})^{\prime \prime }=\bigoplus_{\omega \in {\cal M}%
}B(h_\omega )=\bigoplus_{\omega \in {\cal M}}\pi _\omega ({\cal B}%
)^{\prime \prime }
\]
[KR; 10.3.10]. However we already know that for each $\omega \in {\cal M}$%
, we may consider $h_\omega $ to be a subspace of a possibly larger Hilbert
space $k_\omega $ such that $\pi _\omega \circ \varphi $ is of the form $%
E_\omega \psi _\omega E_\omega $ where $E_\omega $ is the orthogonal
projection of $k_\omega $ onto $h_\omega $, and $\psi _\omega $ is a Jordan $%
*$-morphism from ${\cal A}$ into $B(k_\omega )$ such that the $C^{*}$%
-algebra generated by $\psi _\omega ({\cal A})$ is irreducible, that is $%
\psi _\omega ({\cal A})^{\prime \prime }=B(k_\omega )$ [BR; 2.3.8]. Now
let ${\cal R}$ be the von Neumann algebra $\bigoplus_{\omega \in {\cal %
M}}B(k_\omega ),$ $E$ the projection $\bigoplus_{\omega \in {\cal M}%
}E_\omega $ and $\psi $ the map $\bigoplus_{\omega \in {\cal M}}\psi
_\omega $. Since each $\psi _\omega $ is a Jordan $*$-morphism it can
readily be verified that the same is true of $\psi $. Moreover 
\[
\pi \circ \varphi =\bigoplus_{\omega \in {\cal M}}\pi _\omega \circ
\varphi =\bigoplus_{\omega \in {\cal M}}E_\omega \psi _\omega (\cdot
)E_\omega =E\psi (\cdot )E.
\]
Finally for any $\omega _0\in {\cal M}$, let $F_{\omega _0}$ be the
orthogonal projection of $\bigoplus_{\omega \in {\cal M}}k_\omega $ onto
the subspace cor\-re\-spon\-ding to $k_{\omega _0}$. Clearly the $F_\omega $'s are
ortho\-gonal pro\-jec\-tions be\-lon\-ging to the cen\-tre of 
${\cal R}$ such that $%
\vee _{\omega \in {\cal M}}F_\omega = I$ with $(F_\omega \psi ({\cal A}%
)F_\omega )^{\prime \prime }\equiv \psi _\omega ({\cal A})^{\prime \prime
}=B(k_\omega )\equiv F_\omega {\cal R}F_\omega $ for each $\omega \in 
{\cal M}$. \hfill $\Box $

An analysis of the proof reveals a measure of dependence on some form of
normality for the pure states. For this reason Theorem 5 was stated in terms
of the reduced atomic representation. The following result puts the matter
in context and enables us to restate Theorem 5 in terms of any faithful
representation of ${\cal B}$ in which pure states are normal.

{\bf Proposition 6} \quad Let ${\cal A} \subset B(h)$ be a concrete $%
C^\ast$-algebra with the property that ${\cal P_A} \subset {\cal N_A}$%
. Then there exists a projection $E \in {\cal A}^{\prime}\cap {\cal A}%
^{\prime\prime}$ such that ${\cal A}_E$ affords an ultra-weakly
continuous $\ast$-isomorphic copy of the reduced atomic representation of $%
{\cal A}$.

{\bf Proof} \quad Let $\omega \in {\cal P_A}$ and let $\tilde{\omega}$
be the unique ultra-weakly continuous pure state extension of $\omega $ to
all of ${\cal A}^{\prime \prime }$. We first show that regarding the GNS
constructions corresponding to $\omega $ and $\tilde{\omega}$ respectively,
we have $h_\omega =h_{\tilde{\omega}}$ with $\pi _\omega =\pi _{\tilde{\omega%
}}\mid _{{\cal A}}$. This effectively follows from the first part of the
proof of Lemma 4 combined with the uniqueness in [BR; 2.3.16]. To see this
directly note that $\Vert \pi _\omega (A)\Omega _\omega \Vert ^2=\omega
(A^{*}A)=\tilde{\omega}(A^{*}A)=\Vert \pi _{\tilde{\omega}}(A)\Omega _{%
\tilde{\omega}}\Vert ^2$ for all $A\in {\cal A}$ where $\Omega _\omega $,
\ $\Omega _{\tilde{\omega}}$ are the relevant canonical cyclic vectors. It
is obvious that the set $\{\pi _{\tilde{\omega}}(A)\Omega _{\tilde{\omega}%
}:A\in {\cal A}\}\subset h_{\tilde{\omega}}$ affords an isometric copy of
the dense subset $\{\pi _\omega (A)\Omega _\omega :A\in {\cal A}\}$ of $%
h_\omega $. If we can show that the former is also dense in $h_{\tilde{\omega%
}}$, then surely $h_\omega =h_{\tilde{\omega}}$ in a canonical way, in which
case we are done as regards the first part of the proof. To see this we
first note that $\pi _{\tilde{\omega}}({\cal A}^{\prime \prime })=B(h_{%
\tilde{\omega}})$ by [KR; 7.1.7 \& 10.2.3]. Then surely $\pi _{\tilde{\omega}%
}$ is ultra-weakly continuous [BR; 2.4.23]. Since ${\cal A}$ is
ultra-weakly dense in ${\cal A}^{\prime \prime }$ [BR; 2.4.11], $\pi _{%
\tilde{\omega}}({\cal A})$ must then be ultra-weakly dense in $\pi _{%
\tilde{\omega}}({\cal A}^{\prime \prime })=B(h_{\tilde{\omega}})$ by
continuity, and hence even strongly dense by [BR, 2.4.11]. But then $\pi _{%
\tilde{\omega}}({\cal A})\Omega _{\tilde{\omega}}=\{\pi _{\tilde{\omega}%
}(A)\Omega _{\tilde{\omega}}:A\in {\cal A}\}$ is dense in $h_{\tilde{%
\omega}}=B(h_{\tilde{\omega}})\Omega _{\tilde{\omega}}=\{B\Omega _{\tilde{%
\omega}}:B\in B(h_{\tilde{\omega}})\}$ as required.

If we now apply [BR; 2.4.22 \& 2.4.23] to the kernel of $\pi _{\tilde{\omega}%
}$, the existence of a projection $F_\omega \in {\cal A}^{\prime }\cap 
{\cal A}^{\prime \prime }$ such that $F_\omega {\cal A}^{\prime \prime
}F_\omega =\pi _{\tilde{\omega}}^{-1}(0)$ follows. Thus denoting $F_\omega
^{\perp }$ by $E_\omega $, it follows that $\pi _{\tilde{\omega}}$ maps $%
E_\omega {\cal A}E_\omega \ ({\rm respectively}\ E_\omega {\cal A}%
^{\prime \prime }E_\omega )$ $*$-isomorphically onto $\pi _\omega ({\cal A%
})\ ({\rm respectively}\ \pi _\omega ({\cal A})^{\prime \prime }=B(h_{%
\tilde{\omega}})=\pi _{\tilde{\omega}}({\cal A}^{\prime \prime }))$. If
now $\rho $ is another pure state on ${\cal A}$ and in the same fashion
we obtain $E_\rho \in {\cal A}^{\prime }\cap {\cal A}^{\prime \prime }$
so that $E_\omega E_\rho \neq 0$, then $E_\omega E_\rho {\cal A}E_\omega
E_\rho \subset (E_\omega {\cal A}E_\omega \cap E_\rho {\cal A}E_\rho
), $ and hence $E_\omega E_\rho {\cal A}E_\omega E_\rho $ affords
equivalent subrepresentations of $\pi _{\tilde{\omega}}({\cal A})=\pi
_\omega ({\cal A})$ and $\pi _{\tilde{\rho}}({\cal A})=\pi _\rho (%
{\cal A})$. Thus by [KR; 10.3.4], $\pi _\omega $ and $\pi _\rho $ are
then not disjoint. If now ${\cal M}$ is a maximal family of pure states
on ${\cal A}$ such that the associated irreducible representations are
pairwise inequivalent, then by [KR; 10.3.7] and the above, the projections $%
\{E_\omega :\omega \in {\cal M}\}\subset {\cal A}^{\prime }\cap 
{\cal A}^{\prime \prime }$ are pairwise orthogonal. Finally since each $%
E_\omega {\cal A}E_\omega ,\ \omega \in {\cal M}$, affords a copy of
the irreducible representation $\pi _\omega ({\cal A})$ (with $E_\omega 
{\cal A}^{\prime \prime }E_\omega $ corresponding to $\pi _{\tilde{\omega}%
}({\cal A}^{\prime \prime })=B(h_{\tilde{\omega}}))$, it follows that
with $E=\bigoplus_{\omega \in {\cal M}}E_\omega \in {\cal A}^{\prime
}\cap {\cal A}^{\prime \prime }$, 
\[
E{\cal A}E=(\bigoplus_{\omega \in {\cal M}}E_\omega ){\cal A}%
(\bigoplus_{\omega \in {\cal M}}E_\omega )=\bigoplus_{\omega \in {\cal %
M}}E_\omega {\cal A}E_\omega
\]
affords a copy of the reduced atomic representation $\bigoplus_{\omega \in 
{\cal M}}\pi _\omega ({\cal A})$ with respect to the map $%
\bigoplus_{\omega \in {\cal M}}\pi _{\tilde{\omega}}$, such that $E%
{\cal A}^{\prime \prime }E=\bigoplus_{\omega \in {\cal M}}E_\omega 
{\cal A}^{\prime \prime }E_\omega $ affords a copy of $\bigoplus_{\omega
\in {\cal M}}B(h_{\tilde{\omega}})$ with respect to the same map. \hfill $%
\Box $

\section{A structural characterisation of maps whose adjoints pre\-ser\-ve
ex\-tre\-me points of the dual ball}

Having dealt with linear maps from one $C^\ast$-algebra into another which
preserve pure states on composition, we now turn our attention to those maps
which preserve the extreme points of the unit ball of the dual of the range
space. We shall eventually see that the ``pure state preserving'' maps have
this property. As might be expected this study does however require a number
of not insubstantial lemmas, the first of which is based on [KR; 7.3.2]. The
fundamental idea behind the first cycle of lemmas is to describe ``extremal
functionals'' in terms of pure states. In so doing we are then able to make
use of the results of \S 1 to achieve the stated objective of this section.

{\bf Lemma 7} \quad Let $\rho$ be a norm-one functional on a $C^\ast$%
-algebra ${\cal A}$. If now we identify ${\cal A}$ with its universal
representation and extend $\rho$ to ${\cal A}^{\prime\prime}= {\cal A}%
^{\ast\ast}$, then by [KR; 7.3.2 \& 10.1.2] there exists a partial isometry $%
V$ in ${\cal A}^{\prime\prime}$ and a normal state $\omega$ (that is $%
\omega \in {\cal A}^\ast_+)$ so that $\omega(A) = \rho(VA)$ and $%
\omega(V^\ast A) = \rho(A)$ for all $A \in {\cal A}$. If indeed $\omega$
is a pure state of ${\cal A}$, we may then assume $V$ to be a unitary
element of ${\cal A}$.

{\bf Proof} \quad Let $\pi_\omega$ be the canonical irreducible
representation of ${\cal A}$ engendered by $\omega$. Since now

\setcounter{equation}{0} 
\begin{equation}
\rho (A)=\omega (V^{*}A)=\left\langle \pi _\omega (A)\Omega _\omega ,\pi
_\omega (V)\Omega _\omega \right\rangle \quad \mbox{\rm for all}
  \quad A\in {\cal A}
\end{equation}

with $\Vert \pi _\omega (V)\Omega _\omega \Vert ^2=\left\langle \pi _\omega
(V^{*}V)\Omega _\omega ,\Omega _\omega \right\rangle =\omega (V^{*}V)=\rho
(V)=\omega (I)=1$, there exists a unitary element $U\in \pi _\omega (%
{\cal A})^{\prime \prime }=B(h_\omega )$ with $U\Omega _\omega =\pi
_\omega (V)\Omega _\omega $. On applying [KR; 5.4.5] we conclude that there
exists $H\in {\cal A}$ such that $\pi _\omega (H)=\pi _\omega (H^{*})$,
and that ${\rm exp}(i\pi _\omega (H))$ maps $\Omega _\omega $ onto $\pi
_\omega (V)\Omega _\omega $. Replacing $H$ by $\frac 12(H+H^{*})$ if
necessary, we may assume $H$ to be self-adjoint. If now we select a sequence
of polynomials $p_n$ such that $p_n\rightarrow {\rm exp}(i\cdot )$
uniformly on $[-\Vert H\Vert ,\Vert H\Vert ]$, then

\[
\begin{array}{lll}
\pi _\omega ({\rm exp}(iH)) & = & \pi _\omega (\lim_np_n(H))=\lim_n\pi
_\omega (p_n(H)) \\ 
& = & \lim_np_n(\pi _\omega (H))={\rm exp}(i\pi _\omega (H))
\end{array}
\]
where the convergence is in norm. Thus with $W={\rm exp}(iH),\ \pi
_\omega (W)$ maps $\Omega _\omega $ onto $\pi _\omega (V)\Omega _\omega $.
>From (1) we then have that 
\[
\rho (A)=\left\langle \pi _\omega (A)\Omega _\omega ,\pi _\omega (W)\Omega
_\omega \right\rangle =\omega (W^{*}A)\quad \mbox{\rm for all}
\quad A\in {\cal A}.
\]
But then $\omega (A)=\omega (IA)=\omega (W^{*}WA)=\rho (WA)$ for all $A\in 
{\cal A}.$ \hfill $\Box $

Although the following is bound to be known, we are not aware of an explicit
reference for it.

{\bf Lemma 8} \quad For any $C^\ast$-algebra ${\cal A}$, \ ${\cal %
S_A}$ is a face of ${\cal A}_1$.

{\bf Proof} \quad Let $\omega $ be a state of ${\cal A}$ and $\rho _1$
and $\rho _2$ functionals in ${\cal A}_1$ with $\lambda \rho
_1+(1-\lambda )\rho _2=\omega $ for some $\lambda $ between $0$ and $1$.
Since then 
\[
\begin{array}{lll}
\omega =\ {\rm Re}(\omega ) & = & \frac{1}{2}((\lambda \rho _1+(1-\lambda
)\rho _2)+(\lambda \rho _1+(1-\lambda )\rho _2)^{*}) \\ 
& = & \lambda (\frac{1}{2}(\rho _1+\rho _1^{*}))+(1-\lambda )(\frac{1}{2}(\rho
_2+\rho _2^{*})) \\ 
& = & \lambda {\rm Re\ }\rho _1+(1-\lambda ){\rm Re\ }\rho _2,
\end{array}
\]
it follows that 
\[
\begin{array}{lll}
1=\ \omega (I) & = & \lambda {\rm Re}(\rho _1(I))+(1-\lambda ){\rm Re}%
(\rho _2(I)) \\ 
& \leq & \lambda |\rho _1(I)|+(1-\lambda )|\rho _2(I)| \\ 
& \leq & \lambda \Vert \rho _1\Vert +(1-\lambda )\Vert \rho _2\Vert \\ 
& \leq & \lambda +(1-\lambda )=1.
\end{array}
\]

But this can only be if 
\[
1 = {\rm Re} (\rho_k(I)) = | \rho_k(I) | = \Vert \rho_k \Vert \quad k =
1,2
\]
in which case 
\[
\rho_k(I) = 1 = \Vert \rho_k \Vert \quad k = 1,2.
\]

Thus as required $\rho_1$ and $\rho_2$ are states by [KR; 4.3.2]. \hfill 
$\Box$

{\bf Corollary 9} \quad If $\rho$ is a bounded functional on a $C^\ast$%
-algebra ${\cal A}$ related to a pure state $\omega$ in the manner
described in the hypothesis of Lemma 7, then $\rho$ is an extreme point of $%
{\cal A}^\ast_1.$

{\bf Proof} \quad First of all note that $\Vert \rho \Vert \leq \Vert
\omega \Vert \Vert V \Vert = 1$. Now suppose $\rho_1$ and $\rho_2$ are
elements of ${\cal A}_1$ such that $\rho = \lambda \rho_1 + (1 - \lambda)
\rho_2$. Then surely 
\[
\omega(A) = \lambda \rho_1(V^\ast A) + (1 - \lambda) \rho_2(V^\ast A) \quad 
\mbox{\rm for all} \quad A \in {\cal A}.
\]

Considering Lemma 8 alongside the fact that $\omega$ is pure, we have $%
\omega(A) = \rho_k (V^\ast A)$ for all $A \in {\cal A}$ and $k = 1,2$.
But then $\rho(A) = \rho_k (VV^\ast A)$ for all $A \in {\cal A}$. Finally
since by Lemma 7 we may assume $V$ to be unitary in ${\cal A}$, we
therefore have that $\rho = \rho_1 = \rho_2$ as required. \hfill $\Box$

We now see that in fact the converse of Corollary 9 also holds. This enables
us to use the theory concerning pure states to treat the extreme points of
the unit ball of the dual of a $C^\ast$-algebra.

{\bf Lemma 10} \quad Let $\rho, \omega, {\cal A}$ and $V$ be as in
Lemma 7. Then $\omega$ is pure whenever $\rho$ is an extreme point of $%
{\cal A}^\ast_1.$

{\bf Proof} \quad Assume $\rho$ to be an extreme point of ${\cal A}%
^\ast_1$, and let $\omega_1, \omega_2$ be states on ${\cal A}$ with

\setcounter{equation}{0} 
\begin{equation}
\omega = \lambda \omega_1 + (1 - \lambda) \omega_2 \quad 0 < \lambda < 1.
\end{equation}

If now we identify $\rho ,\omega ,\omega _1$ and $\omega _2$ with their
canonical ultra-weakly continuous extensions to ${\cal A}^{\prime \prime
}={\cal A}^{**}$ [KR; 10.1.1], then surely 
\[
\rho (A)=\omega (V^{*}A)=\lambda \omega _1(V^{*}A)+(1-\lambda )\omega
_2(V^{*}A)\quad \mbox{\rm for all}\quad A\in {\cal A}%
^{\prime \prime }.
\]

Since $\rho$ is an extreme point of ${\cal A}^\ast_1$ with $\Vert
\omega_k (V^\ast \cdot) \Vert \leq \Vert \omega_k \Vert \Vert V^\ast \Vert =
1$ for $k = 1,2,$ we conclude that $\omega_k(V^\ast A) = \rho(A)$ for all $A
\in {\cal A}^{\prime\prime}$ where $k = 1,2.$ But then

\begin{equation}
\omega (A)=\rho (VA)=\omega _k(V^{*}VA)\quad \mbox{\rm for all}
\quad A\in {\cal A}^{\prime \prime },\quad k=1,2.
\end{equation}

Moreover since $0 \leq \lambda \omega_1 \leq \omega$ and $0 \leq (1 -
\lambda) \omega_2 \leq \omega$ by (1), the fact that $I - V^\ast V \geq 0$
implies that 
\[
0 \leq \lambda \omega_1 (I - V^\ast V) \leq \omega (I - V^\ast V) = 1 - \rho
(V) = 1 - \omega(I) = 0.
\]

Similarly $(1 - \lambda) \omega_2 (I - V^\ast V) = 0,$ and hence 
\[
\omega_k (I - V^\ast V) = 0 \quad {\rm for} \quad k = 1,2.
\]

But then for any $A \in {\cal A}^{\prime\prime}$, the fact that $I -
V^\ast V$ is a projection considered alongside [KR; 4,3,1] implies that 
\[
\begin{array}{lll}
0 \leq |\omega_k ((I - V^\ast V) A)|^2 & \leq & \omega_k ((I - V^\ast V) (I
- V^\ast V)^\ast) \omega_k (A^\ast A) \\ 
& = & \omega_k (I - V^\ast V) \omega_k (A^\ast A) = 0 \quad k = 1,2.
\end{array}
\]

This fact combined with (2) now reveals that $\omega = \omega_1 = \omega_2$
as required. \hfill $\Box$

Having achieved the objective of describing extreme points of ${\cal A}%
^\ast_1$ in terms of pure states, we are now able to duplicate the
fundamental lemmas of \S 1 for the more general case.

{\bf Lemma 11} \quad Let ${\cal A}$ be a concrete $C^\ast$-algebra
with ${\cal P_A} \subset N_{{\cal A}}$ (eg. the reduced atomic
representation). Then every extreme point $\rho$ of ${\cal A}^\ast_1$ has
a unique ultra-weakly continuous norm-preserving extension $\tilde{\rho}$ to
all of ${\cal A}^{\prime\prime}$. Moreover $\tilde{\rho}$ is an extreme
point of $({\cal A}^{\prime\prime})^{\ast}_1$.

{\bf Proof} \quad By Lemmas 7 and 10, and Corollary 9 there exists a
unitary $V \in {\cal A}$ so that $\rho(V \cdot) = \omega$ is a pure state
of ${\cal A}$. By (Lemma 3 and) [KR; 10.1.11], $\omega$ has a unique norm
preserving ultra-weakly continuous extension $\tilde{\omega}$ to all of $%
{\cal A}^{\prime\prime}$. Clearly $\tilde{\omega} (V^\ast A) = \tilde{\rho%
}(A)$ for all $A \in {\cal A}^{\prime\prime}$ then defines a norm
preserving extension $\tilde{\rho}$ of $\rho$ to all of ${\cal A}%
^{\prime\prime}$, which is moreover ultra-weakly continuous by the
ultra-weak continuity of $\tilde{\omega}$ combined with [BR; 2.4.2]. But
then $\tilde{\rho}|_{{\cal A}} = \rho$ is ultra-weakly continuous on $%
{\cal A}$, and so the extension $\tilde{\rho}$ must be unique by [KR;
10.1.11]. Finally since $V$ is unitary and since $\tilde{\omega}$ is a pure
state of ${\cal A}^{\prime\prime}$ by Lemma 4, it now follows from
Corollary 9 that $\tilde{\rho}$ is an extreme point of $({\cal A}%
^{\prime\prime})^{\ast}_1$. \hfill $\Box$

{\bf Lemma 12} \quad Let ${\cal A} \subset B(h)$ be a concrete $C^\ast$%
-algebra, $E$ a projection in ${\cal A}$, and let $\eta$ be defined as in
Lemma 1. Then for any ultra-weakly continuous functional $\rho$ on ${\cal %
A}_E$, \ $\rho \circ \eta$ is an ultra-weakly continuous functional on $%
{\cal A}$ with $\Vert \rho \Vert = \Vert \rho \circ \eta \Vert.$
Conversely if $\tilde{\rho}$ is ultra-weakly continuous on ${\cal A}$,
then the restriction of $\tilde{\rho}$ to ${\cal A}$ is ultra-weakly
continuous with respect to ${\cal A}_E$.

{\bf Proof} \quad This is a fairly obvious and easily verifiable
consequence of Lemma 2 considered alongside the fact that each ultra-weakly
continuous functional is a linear combination of normal states (see for
example [KR; 7.4.7]). We therefore forgo the proof. \hfill $\Box$

{\bf Lemma 13} \quad Let ${\cal A}$ be a $C^\ast$-algebra, $E$ a
projection in ${\cal A}$, and $\rho$ an extreme point of $({\cal A}%
^\ast_E)_1$. With $\eta$ defined as in Lemma 1, it then follows that $\rho
\circ \eta$ is an extreme point of ${\cal A}^\ast_1$. Conversely if $%
\tilde{\rho}$ is an extreme point of ${\cal A}^\ast_1$ with $\Vert \tilde{%
\rho} |_{{\cal A}_E} \Vert = 1$, then the restriction of $\tilde{\rho}$
to ${\cal A}_E$ is an extreme point of $({\cal A}^\ast_E)_1$.

{\bf Proof} \quad By Lemmas 7 and 10, and Corollary 9 there exists a
unitary element $V$ of ${\cal A}_E$ such that $\omega = \rho(V \cdot)$ is
a pure state of ${\cal A}_E$. By Lemma 1 $\omega \circ \eta$ is a pure
state of ${\cal A}$. But since $VE = EV = V$ with $V^\ast V = VV^\ast = E$%
, it is clear that $\tilde{V} = V + (I - E)$ is unitary in ${\cal A}$.
Moreover since then 
\[
\omega \circ \eta (A) = \rho(EVAE) = \rho(E\tilde{V} AE) = \rho \circ \eta (%
\tilde{V} A)
\]
for all $A \in {\cal A}$, it is fairly clear from Corollary 9 that $\rho
\circ \eta$ is then an extreme point of ${\cal A}^\ast_1$. Conversely
suppose $\tilde{\rho}$ is an extreme point of ${\cal A}^\ast_1$ with $%
\Vert \tilde{\rho} |_{{\cal A}_E} \Vert = 1$ and assume that ${\cal A}
\subset B(h)$ is universally represented. By Lemma 11 we may identify $%
\tilde{\rho}$ with its unique ultra-weakly continuous extension to ${\cal %
A}^{\prime\prime}$. But then Lemma 12 informs us that the restriction $%
\tilde{\rho} |_{{\cal A}_E}$, which we will henceforth denote by $\rho_0$%
, is an ultra-weakly continuous functional on ${\cal A}^{\prime\prime}_E$%
, with $\Vert \rho_0 \Vert = 1$ by assumption. Applying [KR; 7.3.2], we
conclude that there exists a partial isometry $V \in {\cal A}%
^{\prime\prime}_E$ and a normal state $\omega$ on ${\cal A}%
^{\prime\prime}_E$ so that $\rho_0(VA) = \omega(A)$ and $\rho_0(A) =
\omega(V^\ast A)$ for all $A \in {\cal A}^{\prime\prime}_E$. But then $EV
= VE = V$ with $\tilde{\rho} (VI) = \tilde{\rho}(VE) = \rho_0 (VE) =
\omega(E) = 1$ by [KR; 4.3.2]. Hence again by [KR; 4.3.2], $\tilde{\rho}%
(VA), \ A \in {\cal A}$, defines a state on ${\cal A}$. By Lemma 10
this state is necessarily a pure state of ${\cal A}$, and hence by Lemma
1, $\omega$ is a pure state of ${\cal A}_E$ on restriction to ${\cal A}%
_E$. Finally on considering [KR; 10.1.12 \& 10.1.21] it is clear that we may
assume $V \in {\cal A}^{\prime\prime}_E \subset ({\cal A}%
_E)^{\ast\ast} $ (up to an isometric isomorphism) and hence by Corollary 9,
on restriction to ${\cal A}_E$,\ $\rho_0$ is an extreme point of $(%
{\cal A}^\ast_E)_1$. \hfill $\Box$

The final building block we need to achieve a general characterisation of
maps with ``extreme point preserving'' adjoints, is that of reducing this
question to the case of maps from say $B(h)$ to $B(k)$, where $h$ and $k$
are Hilbert spaces. It seems that we need to take steps to ensure that all
pure states are normal in order to achieve this.

{\bf Lemma 14} \quad Let ${\cal A}$ be a $C^\ast$-algebra and $\rho_1,
\rho_2$ extreme points of ${\cal A}^\ast_1$ for which the associated
(pure) states defined as in Lemma 7 are disjoint, then $\Vert \rho_1 -
\rho_2 \Vert = 2$.

{\bf Proof} \quad Observe that Lemma 10 ensures that the associated
states, say $\omega _1$ and $\omega _2$, are indeed pure. Moreover by Lemma
7 there exist unitaries $U_1,U_2$ in ${\cal A}$ so that $\rho
_i(A)=\omega _i(U_i^{*}A)$ and $\rho _i(U_iA)=\omega _i(A)$ for all $A\in 
{\cal A}$, \ $i=1,2$. Now let $\pi _{\omega _i},\ i=1,2,$ be the
representations engendered by $\omega _i,\ i=1,2$. If now we let $\pi =\pi
_{\omega _1}\oplus \pi _{\omega _2}$ and if indeed $\omega _1$ and $\omega
_2 $ are disjoint, then by [KR; 10.3.3(iii)] with $E^{\prime }$ and $%
F^{\prime } $ the projections of $h_1\oplus h_2$ onto $\{0\}\oplus h_2$ and $%
h_1\oplus \{0\}$ respectively, we surely have 
\[
\pi ({\cal A})=\pi _{\omega _1}({\cal A})\oplus \pi _{\omega _2}(%
{\cal A}).
\]

Hence we may select $V\in {\cal A}$ with $\pi _{\omega _1}(V)=\pi
_{\omega _1}(U_1)$ and $\pi _{\omega _2}(V)=-\pi _{\omega _2}(U_2)$.
Moreover since then $\Vert \pi (V)\Vert =1$, we may in fact select $V$ so
that $\Vert V\Vert =1$, since $\pi ({\cal A}_1)=\pi ({\cal A})_1$.
Thus 
\[
\begin{array}{lll}
2\geq \Vert \rho _1-\rho _2\Vert & \geq & |\rho _1(V)-\rho _2(V)|=|\omega
_1(U_1^{*}V)-\omega _2(U_2^{*}V)| \\ 
& = & |\left\langle \pi _{\omega _1}(U_1^{*}V)\Omega _{\omega _1},\Omega
_{\omega _1}\right\rangle -\left\langle \pi _{\omega _2}(U_2^{*}V)\Omega
_{\omega _2},\Omega _{\omega _2}\right\rangle | \\ 
& = & |\left\langle \pi _{\omega _1}(U_1^{*}U_1)\Omega _{\omega _1},\Omega
_{\omega _1}\right\rangle +\left\langle \pi _{\omega _2}(U_2^{*}U_2)\Omega
_{\omega _2},\Omega _{\omega _2}\right\rangle | \\ 
& = & \left\langle \Omega _{\omega _1},\Omega _{\omega _1}\right\rangle
+\left\langle \Omega _{\omega _2},\Omega _{\omega _2}\right\rangle =2.
\end{array}
\]
\hfill $\Box $

The next Lemma in this cycle is based on an adaptation of a technique of
St\o rmer's [St\o 2; 5.6].

{\bf Lemma 15} \quad Let ${\cal A}$ and ${\cal B}$ be $C^\ast$%
-algebras.

\begin{description}
\item[(a)]  If $\psi :{\cal A}\rightarrow {\cal B}$ is a linear map
for which $\rho \circ \psi \in {\rm ext}({\cal A}_1^{*})$ whenever $%
\rho \in {\rm ext}({\cal B}_1^{*})$, then given any two unitarily
equivalent pure states $\omega _1,\omega _2$ of ${\cal B}$, the pure
states associated with $\omega _1\circ \psi $ and $\omega _1\circ \psi $ by
means of the technique described in Lemma 7, are also unitarily equivalent.

\item[(b)]  If ${\cal B}=B(h)$ and if $\psi :{\cal A}\rightarrow 
{\cal B}$ is a linear map for which $\rho \circ \psi \in {\rm ext}(%
{\cal A}_1^{*})$ whenever $\rho $ is an ultra-weakly continuous element
of ${\rm ext}({\cal B}_1^{*})$, then given any two unitarily
equivalent ultra-weakly continuous pure states $\omega _1$ and $\omega _2$
of ${\cal B}$, the pure states associated with $\omega _1\circ \psi $ and 
$\omega _2\circ \psi $ are unitarily equivalent.
\end{description}

{\bf Proof} \quad The proof of the two cases being very similar, we
content ourselves with proving (a). In fact as regards the proof of (b) as
compared to (a), the only additional piece of information we need for (b) is
to note that the stated condition in (b) is sufficient to ensure that $\Vert
\psi \Vert \leq 1$. To see this note that for any given $\epsilon >0$ and $%
A\in {\cal A}$ we may select $x,y\in h$ with $\Vert x\Vert =\Vert y\Vert
=1$ so that $\Vert \psi (A)\Vert -\epsilon \leq \left| \left\langle \psi
(A)x,y\right\rangle \right| $. From for example Corollary 9 and [KR; 4.6.8]
we may easily deduce that the functional $\rho (T)=\left\langle
Tx,y\right\rangle ,\ T\in {\cal B}$, is an ultra-weakly continuous
element of ${\rm ext}({\cal B}_1^{*})$. But then $\rho \circ \psi \in 
{\rm ext}({\cal A}_1^{*})$, and hence $\Vert \psi (A)\Vert -\epsilon
\leq \left| (\rho \circ \psi )(A)\right| \leq \Vert \rho \circ \psi \Vert
\Vert A\Vert =\Vert A\Vert $. Now for (a) suppose there exists a unitary $%
U\in {\cal B}$ so that $\omega _2(A)=\omega _1(U^{*}AU)$ for all $A\in 
{\cal B}$ where $\omega _1,\omega _2$ are pure states. If now $\pi _1$ is
the canonical representation engendered by $\omega _1$ on say $h_1$ with
corresponding cyclic vector $\Omega $, then surely $\omega _1=\omega _\Omega
\circ \pi _1$ and $\omega _2=\omega _z\circ \pi _1$ where $\omega _\Omega $
and $\omega _z$ are the vector states on $B(h_1)$ corresponding to $\Omega $
and $z=\pi _1(U)\Omega $ respectively. The rest of the proof is basically an
adaptation of part of [St\o 2; 5.6]. Now if $\pi _1(U)\Omega $ was merely a
(modulus one) scalar multiple of $\Omega $, it trivially follows that $%
\omega _\Omega =\omega _z$, and hence in this case we are done since then $%
\omega _1=\omega _z\circ \pi _1=\omega _2$. Thus suppose ${\rm span}%
\{\Omega ,\pi _1(U)\Omega \}=k$ is a two-dimensional subspace of $h_1$, and
select $x\in k$ so that $x\perp \Omega $ with $\Vert x\Vert =1$. Next select 
$\lambda \in {\bf C},\ |\lambda |=1$ so that $\lambda \left\langle \pi
_1(U)\Omega ,\Omega \right\rangle =\left| \left\langle \pi _1(U)\Omega
,\Omega \right\rangle \right| $. Since the vector state induced by $\lambda
\pi _1(U)\Omega $ is identical to $\omega _z$, it follows that we may assume 
$\left\langle \pi _1(U)\Omega ,\Omega \right\rangle =\left| \left\langle \pi
_1(U)\Omega ,\Omega \right\rangle \right| $. By similarly adjusting $x$ if
necessary, we may assume $\left\langle \pi _1(U)\Omega ,x\right\rangle
=\left| \left\langle \pi _1(U)\Omega ,x\right\rangle \right| $. Now let $%
w_1=\Omega ,\ w_2=2^{-\frac 12}(\Omega +x)$ and $w_3=\pi _1(U)\Omega $.
Since $\{\Omega ,x\}$ is an ONB for $k$, it is an easy exercise to show that 
$\Vert w_2\Vert =1$ with $\Vert w_1-w_2\Vert ^2=2-\sqrt{2}<1$. Moreover this
fact together with the foregoing implies that 
\[
w_3=\left\langle w_3,\Omega \right\rangle \Omega +\left\langle
w_3,x\right\rangle x=\left| \left\langle w_3,\Omega \right\rangle \right|
\Omega +\left| \left\langle w_3,x\right\rangle \right| x.
\]

But as $\left| \left\langle w_3,\Omega \right\rangle \right| \leq 1$ and $%
\left| \left\langle w_3,x\right\rangle \right| \leq 1$, we therefore have
that 
\[
\left| \left\langle w_3,\Omega \right\rangle \right| +\left| \left\langle
w_3,x\right\rangle \right| \geq \left| \left\langle w_3,\Omega \right\rangle
\right| ^2+\left| \left\langle w_3,x\right\rangle \right| ^2=\Vert w_3\Vert
^2=1.
\]

Thus since $\left\langle w_3,\Omega \right\rangle \geq 0$ and $\left\langle
w_3,x\right\rangle \geq 0$ by construction, 

\[
\begin{array}{lll}
\Vert w_2-w_3\Vert ^2 & = & \Vert w_2\Vert ^2-2^{-\frac 12}\cdot 2(\left|
\left\langle w_3,\Omega \right\rangle \right| +\left| \left\langle
w_3,x\right\rangle \right|) +\Vert w_3\Vert ^2 \\ 
& = & 2-2^{\frac 12}\left( \left| \left\langle w_3,\Omega \right\rangle \right|
+\left| \left\langle w_3,x\right\rangle \right| \right) \\ 
& \leq & 2-2^{\frac 12}<1.
\end{array}
\]

Let $\nu _i$ be the vector state on ${\cal B}(h_1)$ engendered by $w_i,\
i=1,2,3$. Clearly $\nu _i\circ \pi _1$ is a pure state of ${\cal B}$ for
all $i=1,2,3$ [KR; 10.2.3 \& 10.2.5] with $\nu _1\circ \pi _1=\omega _1$ and 
$\nu _3\circ \pi _1=\omega _2$. Moreover if the pure states associated with $%
\nu _i\circ \pi _1\circ \psi $ and $\nu _{i+1}\circ \pi _1\circ \psi ,\
i=1,2 $, are unitarily equivalent, the same is trivially true of $\omega
_1\circ \psi $ and $\omega _2\circ \psi $. Finally observe that for any $%
i=1,2$, and any $A\in {\cal A}$, we have 
\[
\begin{array}{ll}
\left| \nu _i\circ \pi _1\circ \psi (A)-\nu _{i+1}\circ \pi _1\circ \psi
(A)\right| &  \\ 
\hspace{2cm}=\left| \left\langle \pi _1\circ \psi (A)w_i,w_i\right\rangle
-\left\langle \pi _1\circ \psi (A)w_{i+1},w_{i+1}\right\rangle \right| &  \\ 
\hspace{2cm}=\left| \left\langle \pi _1\circ \psi
(A)(w_i-w_{i+1}),w_i\right\rangle +\left\langle \pi _1\circ \psi
(A)w_{i+1},w_i-w_{i+1}\right\rangle \right| &  \\ 
\hspace{2cm}\leq 2\Vert \pi _1\Vert \Vert \psi \Vert \Vert w_i-w_{i+1}\Vert
\Vert A\Vert &  \\ 
\hspace{2cm}=\left( 2\Vert \psi \Vert (2-2^{\frac 12})^{\frac 12}
   \right) ~\Vert A\Vert\,
. & 
\end{array}
\]

All that remains to be done is to note that since $\psi^\ast$ preserves
extreme points, we necessarily have $\Vert \psi \Vert \leq 1$, and then to
apply Lemma 14. \hfill $\Box$

{\bf Theorem 16} \quad Let $h,k$ be Hilbert spaces.

\begin{enumerate}
\item[(A)]  A continuous linear map $\psi :K(k)\rightarrow K(h)$ has the
property that $\rho \circ \psi \in {\rm ext}({\cal S}_1(k)_1)$
whenever $\rho $ is an extreme point of ${\cal S}%
_1(h)_1$, if and only if $\psi $ is of precisely one of the following forms:

\begin{enumerate}
\item[1)]  There exist injective partial isometries $U:h\rightarrow k$ and $%
V:h\rightarrow k$ such that either $\psi (T)=U^{*}TV$ for all $T\in K(k)$ or 
$\psi (T)=U^{*}c^{*}T^{*}cV$ for all $T\in K(k)$. (Here $c:k\rightarrow k$
is the anti-unitary operator induced by complex conjugation of the scalars.)

\item[2)]  There exists a fixed unit vector $w\in k$ and a surjective
partial isometry $V:k\rightarrow {\cal S}_2(h)$ such that either $\psi
(T)=JV(Tw)$ for all $T\in K(k)$ or $\psi (T)=(JV((T^{*})w))^{*}$ for all $%
T\in K(k)$, where $J$ is the natural injection of ${\cal S}_2(h)$ into $%
K(h).$
\end{enumerate}

\item[(B)]  An ultra-weakly continuous linear map $\psi :B(k)\rightarrow
B(h) $ has the property that $\rho \circ \psi \in {\rm ext}(B(k)_1^{*})$
whenever $\rho $ is an ultra-weakly continuous extreme point of $B(h)_1^{*}$%
, if and only if $\psi $ is of precisely one of the following forms:

\begin{enumerate}
\item[1)]  There exist injective partial isometries $U:h\rightarrow k$ and $%
V:h\rightarrow k$ such that either $\psi (T)=V^{*}TU$ for all $T\in B(k)$ or 
$\psi (T)=V^{*}c^{*}T^{*}cU$ for all $T\in B(k)$. (Here $c:k\rightarrow k$
is the anti-unitary operator induced by complex conjugation of the scalars.)

\item[2)]  There exists a fixed unit vector $w\in k$ and a surjective
partial isometry $V:k\rightarrow {\cal S}_2(h)$ such that either $\psi
(T)=JV(Tw)$ for all $T\in B(k)$ or $\psi (T)=(JV((T^{*})w))^{*}$ for all $%
T\in B(k)$, where $J$ is the natural injection of ${\cal S}_2(h)$ into $%
B(h).$
\end{enumerate}
\end{enumerate}

{\bf Proof} \quad (A)\quad Let $(e_\lambda)_\Gamma$ be a fixed
or\-tho\-nor\-mal
ba\-sis for $h$. Since 

\[{\rm ext}\,\left(K(h)_1^{*}\right) =
{\rm ext}\,{\cal S}(h)_1 = \{u\otimes v\,:u,v\; 
\mbox{\rm unit vectors in}\; h\},\] 

we start the investigation by looking at the images 
$\Phi ^{*}(e_\lambda\bigotimes e_\mu)$, $\lambda,\mu\in\Gamma$. 
Let us state a sublemma providing us with a criterion that
will be repeatedly applied in what follows (its proof is a
straightforward exercise):

\begin{description}
\item[Sublemma]  Let $\lambda_1,\lambda_2,\mu_1,\mu_2\in\Gamma$ 
with $\lambda_1\neq\mu_1,\lambda_2\neq\mu_2$, 
and $u,v,w,z\in k$\ be unit vectors. 
If $\psi ^{*}(e_{\lambda_1}\otimes e_{\lambda_2})=u\otimes v$
and $\psi ^{*}(e_{\lambda_1}\otimes e_{\mu_2})=w\otimes z$, then we either have 
that $u\Vert w$\ and $v\perp z$, or $u\perp w$\ and $v\Vert z$. 
Similarly, if $\psi^{*}(e_{\lambda_1}\otimes e_{\lambda_2})=u\otimes v$
and $\psi ^{*}(e_{\mu_1}\otimes e_{\lambda_2})=w\otimes z$, 
then either $u\Vert w$\ and $v\perp z$, or $u\perp w$\ and $v\Vert z$.
\end{description}

On fixing $\lambda_0\in\Gamma$ and applying the sublemma to the sets 
$(\psi^{*}(e_{\lambda_0}\otimes e_\lambda))_\lambda$ 
and $(\psi ^{*}(e_\lambda\otimes e_{\lambda_0}))_\lambda$, we deduce that
there are four possibilities for their values:

(1$_a$) $\psi ^{*}(e_{\lambda_0}\otimes e_\lambda)=u_{\lambda_0}\otimes v_\lambda$ 
and $\psi ^{*}(e_\lambda\otimes e_{\lambda_0})=u_\lambda\otimes v_{\lambda_0}$

(1$_b$) $\psi ^{*}(e_{\lambda_0}\otimes e_\lambda)=u_\lambda\otimes v_{\lambda_0}$ 
and $\psi ^{*}(e_\lambda\otimes e_{\lambda_0})=u_{\lambda_0}\otimes v_\lambda$

(2$_a$) $\psi ^{*}(e_{\lambda_0}\otimes e_\lambda)=u_\lambda\otimes v_{\lambda_0}$ 
and $\psi ^{*}(e_\lambda\otimes e_{\lambda_0})=w_\lambda\otimes v_{\lambda_0}$

(2$_b$) $\psi ^{*}(e_{\lambda_0}\otimes e_\lambda)=u_{\lambda_0}\otimes v_\lambda$ 
and $\psi ^{*}(e_\lambda\otimes e_{\lambda_0})=u_{\lambda_0}\otimes w_\lambda$

\noindent where all $(u_\lambda),(v_\lambda),(w_\lambda)$ 
are orthonormal systems in $k$ (and $%
w_{\lambda_0}=u_{\lambda_0}$ or $v_{\lambda_0}$, depending on the case).

In the following we will assume $(e_\lambda)_\Gamma$ to be countable and
let $\lambda_0=1$. The reason we may do this is that in each case it is
enough to establish the action of $\psi^{*}$ in terms of arbitrary
countable subsets of $(e_\lambda)_\Gamma$ containing $e_{\lambda_0}$,
in order to establish the action of $\psi^{*}$ in terms of all of
$(e_\lambda)_\Gamma$.

Let us look at the case (1$_a$) first. If we assume that
$\psi ^{*}(e_2\otimes e_2)=y\otimes v_1$ with $y\perp u_2$ (one of the
two possibilities allowed by the sublemma), then on comparing
$\psi ^{*}(e_2\otimes e_3)$ to $\psi ^{*}(e_2\otimes e_i)$ ($i=1,2$)
and applying the sublemma, we have that $\psi ^{*}(e_2\otimes e_3)
=z\otimes v_1$ for some unit vector $z$ such that $z\perp y$. But then
the sublemma applied to $\psi^{*}(e_1\otimes e_2)$ 
and $\psi ^{*}(e_2\otimes e_2)$ would give $u_1\Vert y$ 
(since $v_1\perp v_2$ by assumption). Applying it to 
$\psi^{*}(e_1\otimes e_3)$ and $\psi ^{*}(e_2\otimes e_3)$ gives in turn
$u_1\Vert \ z$ (since $v_1\perp v_3$). Hence, we would get $y\parallel z$, a
contradiction. Thus we must have $\psi ^{*}(e_2\otimes e_2)=u_2\otimes z$
where $z$ is a unit vector with $z\perp v_1$. Inductively applying the
sublemma we have

\[\psi ^{*}(e_2\otimes e_j)=u_2\otimes z_j
   \quad\mbox{\rm for all}\quad j\in{\bf N}\]

where $(z_j)$ is an ONS with $z_1=v_1$. More generally we may verify that
for each $m\in{\bf N}$,

\[\psi ^{*}(e_m\otimes e_j)=u_m\otimes z^{(m)}_j\quad 
       \mbox{\rm for all}\quad j\in{\bf N}\]

where $(z^{(m)}_j)_j$ is an ONS with $z^{(m)}_1=v_1$. Thus by symmetry it
follows that there is only the following way to satisfy
the sublemma in this case: we have $\psi ^{*}(e_i\otimes e_j)={\varepsilon }%
_{ij}u_i\otimes v_j$ for all $i,j$, where the $\varepsilon _{ij}$ are
complex numbers of modulus one (we set $\varepsilon _{ij}:=1$ whenever $i$
or $j$ are 1). Define $U,V$ (injective) partial isometries $h\rightarrow k$
by $Ue_i:=u_i$ and $Ve_i:=v_i$ for all $i$. We then have

\[
V^{*}\psi ^{*}(e_i\otimes e_j)U=\varepsilon _{ij}V^{*}(u_i\otimes
v_j)U=\varepsilon _{ij}e_i\otimes e_j\quad\mbox{\rm for all}
     \quad j\,.
\]

Call $\Psi :{\cal S}_1(h)\to {\cal S}_1(h),\,S\mapsto V^{*}\psi
^{*}(S)U$. Then $\left\| \Psi \right\| =1$ and $\Psi $ acts on ${\cal S}%
_1(h)$ as a Schur multiplier with matrix $(\varepsilon _{ij})$ (with respect
to the basis $(e_i)$, of course). It is not difficult to prove that this
contradicts $\left\| \Psi \right\| \leq 1$ (even in the case of
two-dimensional $k$!) except when all the $\varepsilon _{ij}$ are the same
unimodular number, i.e. equal to $\varepsilon _{11}=1$. Accordingly, we can
now write that $\psi ^{*}(e_i\otimes e_j)=V(e_i\otimes e_j)U^{*}$ for all $%
i,j$ and so

\[
\psi ^{*}(S)=VSU^{*}\quad \mbox{\rm for all}\quad S\in S_1(h)\,.
\]
Finally, if $T\in K(h)$ and $S\in {\cal S}_1(k)$ are arbitrary, then

\[
{\rm tr}(S\psi (T))={\rm tr}(\psi ^{*}(S)T)={\rm tr}(VSU^{*}T)=%
{\rm tr}(SU^{*}TV)
\]
and so

\[
\psi (T)=U^{*}TV\quad\mbox{\rm for all}\quad T\in K(k)\,.
\]

Let us now deal with the case (2$_a$). Suppose that $\psi ^{*}(e_2\otimes
e_2)=w_2\otimes y$ where $y\in k$ is some unit vector orthogonal to $v_1$.
On comparing $\psi ^{*}(e_2\otimes e_3)$ to $\psi ^{*}(e_2\otimes e_1)$
and $\psi ^{*}(e_2\otimes e_2)$ and applying the sublemma, we conclude that 
$\psi ^{*}(e_2\otimes e_3)=w_2\otimes z$ with $z\in k$ a unit vector
orthogonal to both $y$ and $v_1$. Then the sublemma
applied to $\psi ^{*}(e_1\otimes e_2)$ and $\psi ^{*}(e_2\otimes e_2)$ gives
(since $y\perp v_1$) $u_2\parallel w_2$. Keeping this in mind and applying
the sublemma to $\psi ^{*}(e_1\otimes e_2)$ and $\psi ^{*}(e_1\otimes e_3)$
gives $u_3\perp w_2$. Hence, after looking at $\psi ^{*}(e_1\otimes e_3)$
and $\psi ^{*}(e_2\otimes e_3)$ we would get $z\parallel v_1$, a
contradiction with the above. Thus we must have $\psi ^{*}(e_2\otimes e_2)
=z\otimes v_1$ where $z\in k$ is a unit vector orthogonal to $w_2$.
Continuing inductively it follows that

\[\psi ^{*}(e_2\otimes e_j)=z_j\otimes v_1\quad
     \mbox{\rm for all}\quad j\in{\bf N},\]

and more generally that

\[\psi ^{*}(e_m\otimes e_j)=z^{(m)}_j\otimes v_1\quad
    \mbox{\rm for all}\quad j\in{\bf N}\]

for each fixed $m\in{\bf N}$, where $(z^{(m)}_j)_j$ is an ONS 
with $z^{(m)}_1=w_2$. By symmetry it follows that there is only the following
way to satisfy the sublemma in this case: we must have $\psi ^{*}(e_i\otimes
e_j)=u_{ij}\otimes v_1$ for all $i,j$ (we renamed the $u_i$ to $u_{1i}$ and
the $w_i$ to $u_{i1}$), where the vectors in the ``matrix'' $(u_{ij})$ form
orthonormal systems along the rows and columns. Now, since $\left\| \psi
^{*}\right\| \leq 1$, we have

\begin{eqnarray*}
{\left\| \sum_{i=1}^n\alpha _i(\sum_{j=1}^n\beta _ju_{ij})\right\| }_k
&=&{\left\| \sum_{i,j=1}^n\,\alpha _i\beta _j(u_{ij}\otimes v_1)
   \right\|}_{{\cal S}_1(k)}\\
&=&\left\| {\psi ^{*}}\left[ {(\sum_{i=1}^n
   \overline{\alpha }_ie_i)}\otimes {(\sum_{i=1}^n\beta _ie_i)}
    \right] \right\|\\
&\leq& (\sum_{i=1}^n\left| \alpha _i\right| ^2)^{1/2}(\sum_{j=1}^n\left|
     \beta _j\right| ^2)^{1/2}
\end{eqnarray*}

for all $(\alpha _i),(\beta _j)\in {\bf C}^n$ and all $n$. Fixing $(\beta
_j) $, this implies that the sequence $(\sum_j\beta _ju_{ij})_{i=1}^n$
must be orthogonal, which in turn --- since $(\beta _j)$ is arbitrary ---
forces the family $(u_{ij})$ to be orthogonal (not only row- and
columnwise!). Let $V:k\to S_2(h)$ be the surjective partial isometry
such that $V^{*}(e_i\otimes e_j)=u_{ij}$ for all $i,j$. Then, if $J$ is the
natural injection ${\cal S}_2(h)\to K(h)$, we have

\[
\psi (T)=JV\left( Tv_1\right) \quad \mbox{\rm for all}\quad T\in K(k)\,.
\]

It should now be clear how to proceed in the analysis of the remaining
cases. To finish the proof of (A), one only needs to check that the
conditions on $\psi$ given in the statement are sufficient to ensure 
$\psi^{*}({\rm ext}\,K(h)_1^{*})\subset {\rm ext}\,K(k)_1^{*}$.
Recall first that ${\rm ext}\,\left( K(h)_1^{*}\right) =
\{u\otimes v\,:u,v$ unit vectors in $h\}$ 
(the same being true for $h$ replaced
by $k$). Suppose then that 1) holds and $\psi$ is an operator
$K(k)\to K(h)$ such that $\psi(T)=U^{*}TV$ for all $T\in K(k)$ and for some
fixed injective partial isometries $U:h\to k$ and $V:h\to k$. This implies
that 

\[\psi^{*}(u\otimes v) = V(u\otimes v)U^{*} = Uu\otimes Vv
\in{\rm ext}\,K(k)_1^{*}\]

for every pair of unit vectors $u,v\in h$. If
$\psi$ were of the form $\psi (T) = U^\ast c^\ast T^\ast cV$ the argument 
would be
entirely similar. On the other hand, if 2) is satisfied and $\psi(T)=
JV(Tw)$ for some fixed unit vector $w\in k$, some fixed surjective partial 
isometry $V:k\to {\cal S}_2(h)$ and all $T\in K(k)$ (and $J$
the natural injection ${\cal S}_2(h)\to K(h)$), then

\[\psi^{*}(u\otimes v)=V^{*}(u\otimes v)\otimes w\in {\rm ext}\,K(k)_1^{*}\]

for all unit vectors $u,v\in h$. Again, if instead $\psi$ were of the form
$\psi(T)=(JV((T^{*})w))^{*}$ the argument would be similar.

(B)\quad All we need to note is that on $B(h)$ and $B(k)$ the ultra-weak and
weak$^{*}$ topologies coincide. Hence since $K(h)^{**}=B(h)$ and $%
K(k)^{**}=B(k)$ and since by hypothesis $\psi $ is weak$^{*}$-continuous, it
follows that on restriction to $K(h)^{*},\ \psi ^{*}$ is a well defined map
from $K(h)^{*}$ into $K(k)^{*}.$ Moreover this restriction maps ${\rm ext}%
(K(h)_1^{*})$ into ${\rm ext}(K(k)_1^{*})$ if and only if $\psi ^{*}$
maps the weak$^{*}$-continuous elements of ${\rm ext}(B(h)_1^{*})$ into
(the weak$^{*}$-continuous elements of) ${\rm ext}(B(k)_1^{*})$. The
result now follows by applying a similar argument as was used in proving
(A). \hfill $\Box $

{\bf Lemma 17} \quad Let ${\cal A}$ be a $C^{*}$-algebra and $\psi :%
{\cal A}\rightarrow B(h)$ a linear map. Then $\psi ^{*}$ maps the
ultra-weakly continuous extreme points of $B(h)_1^{*}$ into ${\rm ext}(%
{\cal A}_1^{*})$ if and only if there exists a Hilbert space $k$ such
that $\psi =\tilde{\psi}\circ \pi $, where $\pi ({\cal A})\subset B(k)$
is an irreducible representation of ${\cal A}$ on $k$ and $\tilde{\psi}%
:B(k)\rightarrow B(h)$ is an ultra-weakly continuous linear map with the
property that $\rho \circ \tilde{\psi}\in {\rm ext}(B(k)_1^{*})$ whenever 
$\rho $ is an ultra-weakly continuous element of ${\rm ext}(B(h)_1^{*})$.

{\bf Proof} \quad First assume $\psi $ to be of the form $\psi =\tilde{%
\psi}\circ \pi $. By the hypothesis all we then really need to check is that 
$\rho \circ \pi \in {\rm ext}({\cal A}_1^{*})$ whenever $\rho $ is an
ultra-weakly continuous extreme point of $B(k)$. By [BR; 2.4.6] and the
extremality of $\rho $, \ $\rho $ is of the form $\rho (A)=\left\langle
Ax,y\right\rangle $ for some $x,y\in k$ with $\Vert x\Vert =\Vert y\Vert =1$%
. As in the proof of Lemma 7 we may now select a unitary $U\in {\cal A}$
with $\pi (U)^{*}y=x$. Then 
\[
\rho \circ \pi (UA)=\left\langle \pi (A)x,\pi (U)^{*}y\right\rangle
=\left\langle \pi (A)x,x\right\rangle \quad A\in {\cal A},
\]
defines a pure state on ${\cal A}$ by [KR; 10.2.5], and hence $\rho \circ
\pi $ is an extreme point of ${\cal A}_1^{*}$ by Corollary 9.

Conversely assume that $\psi ^{*}$ maps the ultra-weakly continuous elements
of ${\rm ext}(B(k)_1^{*})$ into ${\rm ext}({\cal A}_1^{*})$. Now
let $\omega _0$ be a fixed vector state of $B(h)$. Since $\omega _0$ is pure
[KR; 4.6.68], $\omega _0\circ \psi \in {\rm ext}({\cal A}_1^{*})$ by
hypothesis. Thus by Lemmas 7 and 10, and Corollary 9 there exists a unitary $%
U\in {\cal A}$ so that $\nu (A)=\omega _0\circ \psi (UA),\ A\in {\cal A%
}$, defines a pure state on ${\cal A}$. Now let $\pi $ be the irreducible
representation of ${\cal A}$ on some Hilbert space $h_\nu =k,$ engendered
by the GNS process applied to $\nu $ [KR; 10.2.3]. If $\Omega $ is the
canonical cyclic vector in $k$ corresponding to $\nu $, then surely

\setcounter{equation}{0} 
\begin{equation}
(\omega _0\circ \psi )(A)=\left\langle \pi (A)\Omega ,\pi (U)\Omega
\right\rangle \quad \mbox{\rm for every}\quad A\in {\cal A}.
\end{equation}

Now let $\omega $ be any other ultra-weakly continuous pure state of $B(h)$.
By [KR; 7.1.12] and the extremality of $\omega $, \ $\omega $ is precisely a
vector state of $B(h)$ and hence unitarily equivalent to $\omega _0$. But
then by Lemma 15(b) the pure state associated with $\omega \circ \psi $ is
unitarily equivalent to $\nu $. Considering Lemma 7, this effectively means
that there exist unitaries $V$ and $W$ in ${\cal A}$ so that $\omega
\circ \psi (VWAW^{*})=\nu (A)$ for all $A\in {\cal A}$. Consequently

\begin{equation}
\omega \circ \psi (A)=\nu (W^{*}V^{*}AW)=\left\langle \pi (A)\pi (W)\Omega
,\pi (VW)\Omega \right\rangle
\end{equation}

for all $A \in {\cal A}$. If now $\pi(A) = 0$, it is clear from the above
that then $\omega \circ \psi(A) = 0$ for all vector states of $B(h)$
(ultra-weakly continuous elements of ${\rm ext}(B(h)^\ast_1)) $. Since by
the polarization identity the vector states of $B(h)$ separate the points of 
$B(h)$, it follows that $\psi(A) = 0$ whenever $\pi(A) = 0,$ i.e. $%
\psi^{-1}(0) \supset \pi^{-1}(0)$. Thus $\psi$ induces a well defined linear
map $\tilde{\psi}$ from ${\cal A} / \pi^{-1}(0)$ into $B(h)$. Since $\pi$
effectively identifies ${\cal A} / \pi^{-1}(0)$ with $\pi({\cal A})$,
we may assume $\tilde{\psi}$ to be acting from $\pi({\cal A})$ into $B(h)$%
, in which case $\psi = \tilde{\psi} \circ \pi$ by construction. Moreover as
was seen in for example (2) above, for any vector state (ultra-weakly
continuous pure state) $\omega$ of $B(h), \ \omega \circ \tilde{\psi}$ is
ultra-weakly continuous on $\pi({\cal A})$. Hence by [KR; 7.1.12], $%
\omega \circ \tilde{\psi}$ is ultra-weakly continuous for every $\omega \in 
{\cal N}(B(h))$. Now if $\rho$ is ultra-weakly continuous, then so is $%
\rho^\ast$ [BR; 2.4.2], and hence by this fact and [KR; 7.4.7], each
ultra-weakly continuous functional may be written as a linear combination of
at most four normal states. Clearly then $\tilde{\psi}^\ast$ maps
ultra-weakly continuous functionals onto ultra-weakly continuous
functionals. We conclude that $\tilde{\psi}$ must be ultra-weakly
continuous. Since $\pi({\cal A})$ is irreducible $(\pi({\cal A}%
)^{\prime\prime}= B(k))$ it follows that $\pi({\cal A})$ is ultra-weakly
dense in $B(k)$ [BR; 2.4.15]. Thus the result follows on noting that $\tilde{%
\psi}$ has a unique ultra-weakly continuous extension to all of $B(k)$ [KR;
10.1.10], which we may identify with $\tilde{\psi}$ itself. \hfill $\Box$

{\bf Lemma 18} \quad Let ${\cal A}$ and ${\cal B}$ be $C^\ast$%
-algebras and $\psi : {\cal A} \rightarrow {\cal B}$ a linear map.
Then $\psi^\ast$ maps ${\rm ext}({\cal B}^\ast_1)$ into ${\rm ext}(%
{\cal A}^\ast_1)$ if and only if for every irreducible representation $%
\pi({\cal B}) \subset B(h)$ of ${\cal B}, \ (\pi \circ \psi)^\ast$
maps the ultra-weakly continuous elements of ${\rm ext}(B(h)^\ast_1)$
into ${\rm ext}({\cal A}^\ast_1)$.

{\bf Proof} \quad Suppose that for every irreducible representation $\pi $
of ${\cal B}$, \ $\pi \circ \psi $ satisfies the relevant condition
stated above. Given any extreme point $\rho $ of ${\cal B}_1^{*}$, Lemmas
10 and 7 imply that for some unitary $V\in {\cal B},\ \omega (A)=\rho
(VA)(A\in {\cal A})$ defines a pure state of ${\cal B}$. If now $(\pi
_\omega ,h_\omega ,\Omega _\omega )$ is the canonical irreducible [KR;
10.2.3] GNS representation engendered by $\omega $, then surely 
\[
\rho (A)=\omega (V^{*}A)=\left\langle \pi _\omega (A)\Omega _\omega ,\pi
_\omega (V)\Omega _\omega \right\rangle \quad A\in {\cal A}.
\]

Thus $\rho \circ \psi $ is of the form $\rho _0 \circ \pi _\omega \circ \psi $
where $\rho _0$ is defined by $A \rightarrow \left\langle A \Omega _\omega
, \pi _\omega(V)\Omega _\omega \right\rangle$, $A \in B(h_\omega )$. Since
now $\rho _0$ is clear\-ly an ul\-tra-wea\-kly con\-ti\-nuous e\-le\-ment 
of ${\rm ext}(B(h_\omega )_1^{*})$ 
(see for example Corollary 9 and [KR; 4.6.68]), the
hypothesis ensures that $\rho \circ \psi =\rho _0\circ \pi _\omega \circ
\psi $ belongs to ${\rm ext}({\cal B}_1^{*})$.

Conversely suppose that for some irreducible representation $\pi ({\cal B}%
)$ of ${\cal B}$ on say $h$ there exists an ultra-weakly continuous
extreme point $\rho $ of $B(h)_1^{*}$ such that $\rho \circ \pi \circ \psi $
does not belong to ${\rm ext}({\cal A}_1^{*})$. The lemma then follows
on verifying that $\rho \circ \pi \in {\rm ext}({\cal A}_1^{*})$. To
see this note that the extremality of $\rho $ alongside [BR; 2.4.6] ensures
that $\rho $ is of the form 
\[
\rho (T)=\left\langle Tx,y\right\rangle \quad T\in B(h)
\]
for some unit vectors $x,y\in h$. Denoting the vector state corresponding to 
$x$ by $\omega _x$, it follows that $\omega _x\circ \pi $ is a pure state of 
${\cal B}$ [KR; 10.2.5]. Finally arguing as in the proof of Lemma 7, we
may select a unitary $V\in {\cal B}$ so that $\pi (V^{*})x=y$. But then 
\[
\begin{array}{lll}
\omega _x\circ \pi (VA) & = & \left\langle \pi (A)x,\pi (V^{*})x\right\rangle
\\ 
& = & \left\langle \pi (A)x,y\right\rangle =\rho \circ \pi (A)\qquad A\in 
{\cal B}.
\end{array}
\]
Thus $\rho \circ \pi \in {\rm ext}({\cal B}_1^{*})$ by Corollary 9. %
\hfill $\Box $

Finally with all the groundwork done, we are now ready to verify the desired
characterisation. A slight drawback regarding this characterisation is the
atomistic description given to the so-called ``degenerate'' part of such maps.
A more global description of such maps would have been desirable, but may
however not be possible. What is immediately obvious is the recognisable
similarity between the commutative case and the non-degenerate part in the
general case.

{\bf Theorem 19} \quad Let ${\cal A}$ and ${\cal B}$ be $C^{*}$%
-algebras, and $\psi :{\cal A}\rightarrow {\cal B}$ a linear operator.
Then the following are equivalent:

\begin{enumerate}
\item[(a)]  $\rho \circ \psi $ is an extreme point of ${\cal A}_1^{*}$
whenever $\rho $ is an extreme point of ${\cal B}_1^{*}$.

\item[(b)]  For any irreducible representation $\pi $ on say $h$, there
exists a Hilbert space $k$ such that $\pi \circ \psi $ is of precisely one
of the following forms:

\begin{enumerate}
\item[1)]  There exist injective partial isometries $U:h\rightarrow k$ and $%
V:h\rightarrow k$ and a $*$-(anti)morphism $\alpha $ from ${\cal A}$ into 
$B(k)$ with irreducible range such that 
\[
\pi \circ \psi (A)=V^{*}\alpha (A)U\quad \mbox{\rm for all}
  \quad A\in {\cal A}.
\]

\item[2)]  There exists an irreducible representation $\alpha $ of ${\cal %
A}$ on $k$, a fixed unit vector $\omega \in k$, and a surjective partial
isometry $V:k\rightarrow {\cal S}_2(h)$ such that either

\begin{enumerate}
\item[i)]  {$\pi \circ \psi (A)=JV(\alpha (A)w)$ for all $A\in {\cal A}$}
\end{enumerate}

or

\begin{enumerate}
\item[ii)]  $\pi \circ \psi (A)=(JV(\alpha (A)^{*}w))^{*}$ for all $A\in 
{\cal A}$.
\end{enumerate}

Here $J$ is the natural injection of ${\cal S}_2(h)$ into $B(h)$.
\end{enumerate}

\item[(c)]  With ${\cal B}$ in its reduced atomic representation there
exists a projection $E\in {\cal B}\cap {\cal B}^{\prime }$ such that $%
\psi $ decomposes into a degenerate part 
\[
\psi _{I-E}:{\cal A}\rightarrow {\cal B}_{I-E}:A\rightarrow (I-E)\psi
(A)(I-E)
\]
and a non-degenerate part 
\[
\psi _E:{\cal A}\rightarrow {\cal B}_E:A\rightarrow E\psi (A)E
\]
each with the following structure:

\begin{enumerate}
\item[1)]  For every irreducible representation $\pi _0({\cal B}%
_{I-E})\subset B(h_0)$ of ${\cal B}_{I-E}$ there exists an irreducible
representation $\alpha _0({\cal A})\subset B(k_0)$ of ${\cal A}$, a
unit vector $\omega _0\in k_0$, and an embedding $V_0$ of $k_0$ into $B(h_0)$
such that $V_0((k_0)_1)$ is ultra-weakly dense in $B(h_0)_1$ with either 
\[
\pi _0\circ \psi _{I-E}(A)=V_0(\alpha _0(A)w_0)\quad \mbox{\rm for all}
\quad A\in {\cal A}
\]
or 
\[
\pi _0\circ \psi _{I-E}(A)=(V_0(\alpha _0(A)^{*}w_0))^{*}\quad 
\mbox{\rm for all}\quad A\in {\cal A}.
\]

\item[2)]  There exists a von Neumann algebra ${\cal R}$ acting on some
Hilbert space $h$, a partial isometry $W\in {\cal R}$ with initial
projection $E_1$ and final projection $E_2$, and a Jordan $\ast$-morphism $%
\varphi :{\cal A}\rightarrow {\cal R}$ such that {\sl up to $\ast$%
-isomorphic equivalence}, ${\cal B}_E^{\prime \prime }$ appears as $%
{\cal R}_{E_1}$ (or alternatively ${\cal R}_{E_2}$) with $\psi _E$ $*$%
-isomorphically equivalent to the mapping 
\[
{\cal A}\rightarrow {\cal R}_{E_1}:A\rightarrow E_1\psi (A)WE_1
\]
(or alternatively 
\[
A\rightarrow {\cal R}_{E_2}:A\rightarrow E_2W\varphi (A)E_2).
\]
In addition $\varphi ({\cal A})$ has the density property that for some
set $(F_\nu )$ of mutually orthogonal projections in ${\cal R}\cap 
{\cal R}^{\prime }$ with $\sum F_\nu =I$, we have $(F_\nu \varphi (%
{\cal A})F_\nu )^{\prime \prime }={\cal R}_{F_\nu }$ for each $\nu $.
\end{enumerate}
\end{enumerate}

{\bf Proof} \quad The equivalence of (a) and (b) is an immediate
consequence of Theorem 16 (B) and Lemmas 17 and 18. We therefore need only
verify that (a) follows from (c), and that (c) follows from (b).

{\bf (c) $\Longrightarrow $ (a)} \quad Suppose that $\psi $ is of the
form described in (c). Given any $\rho \in {\rm ext}({\cal B}_1^{*})$,
apply Lemmas 7 and 10 to obtain a unitary $V\in {\cal B}$ such that $%
A\rightarrow \rho (VA),\ A\in {\cal B}$, defines a pure state of $%
{\cal B}$. First of all notice that by Lemma 11 we may identify $\rho $
(and $\rho (V\cdot ))$ with its unique extension to ${\cal B}^{\prime
\prime }$. Thus with $E$ as in the hypothesis, an application of [KR;
4.3.14] reveals that $\rho (VE)\in \{0,1\}$ (since $E^2=E)$. Hence again by
[KR; 4.3.14] either 
\[
\rho (A)=\rho (V(V^{*}A)\rho (VE)=\rho (V((V^{*}A)E))=\rho (EA)
\]
for all $A\in {\cal B}$ when $\rho (VE)=1$, or similarly 
\[
\rho (A)=\rho ((I-E)A)\quad \mbox{\rm for all}\quad A\in 
{\cal B}
\]
if indeed $\rho (VE)=0$ (that is $\rho (V(I-E))=1).$ Now if $\rho (V(I-E))=1$%
, then $\rho $ effectively annihilates ${\cal B}_E$, and by Lemma 13
defines an extreme point of $({\cal B}_{I-E})_1^{*}$. Now given any
extreme point $\tilde{\rho}$ of $({\cal B}_{I-E})_1^{*}$, related to some
pure state $\tilde{\omega}$ in the manner described in Lemmas 7 and 10, and
Corollary 9, it follows that $\tilde{\rho}$ is of the form 
\[
\tilde{\rho}(A)=\left\langle \pi _0(A)x,y\right\rangle \quad\mbox{\rm for all}
\quad A\in {\cal A}
\]
where $\pi _0({\cal B}_{I-E})\subset B(h_0)$ is the canonical irreducible
representation engendered by $\tilde{\omega}$, and $x$ and $y$ suitable unit
vectors. Since the functional $\omega _{x,y}(A)=\left\langle
Ax,y\right\rangle ,\ A\in B(h_0)$, is ultra-weakly continuous on $B(h_0)$,
it follows from the hypothesis that $\omega _{x,y}$ assumes its norm on $%
V_0((k_0)_1)$, and hence that $\omega _{x,y}\circ V_0$ is a norm one
functional of $k_0$. Similarly the ultra-weak continuity of $\omega
_{x,y}^{*}$ ensures that $\omega _{x,y}^{*}\circ V_0$ is also a norm-one
functional of $k_0$. The two cases being similar we now assume that 
\[
\pi _0\circ \psi _{I-E}(A)=(V_0(A^{*}w_0))^{*}\quad \mbox{\rm for all}
\quad A\in {\cal A}.
\]
Now since $\omega _{x,y}^{*}\circ V_0$ has norm one, there exists a unit
vector $z\in k_0$ so that 
\[
(\omega _{x,y}^{*}\circ V_0)(p)=\left\langle p,z\right\rangle \quad {\rm %
for all}\quad p\in k_0.
\]
But then 
\[
\begin{array}{lll}
\tilde{\rho}\circ \psi _{I-E}(A) & = & (\omega _{x,y}\circ \pi _0\circ \psi
_{I-E})(A) \\ 
& = & \omega _{x,y}((V_0(\alpha _0(A^{*})w_0))^{*}) \\ 
& = & \overline{\omega _{x,y}^{*}(V_0(\alpha _0(A)^{*}w_0))} \\ 
& = & \overline{\left\langle \alpha _0(A)^{*}w_0,z\right\rangle } \\ 
& = & \left\langle \alpha _0(A)z,w_0\right\rangle \quad \mbox{\rm for all}
\quad A\in {\cal A}.
\end{array}
\]

If now as in the proof of Lemma 7 we select a unitary $U\in {\cal A}$ so
that $\pi _0(U^{*})w_0=z$, then $\tilde{\rho}\circ \psi
_{I-E}(UA)=\left\langle \alpha _0(U)\alpha _0(A)z,\omega _0\right\rangle
=\left\langle \alpha _0(A)z,z\right\rangle ,\ A\in {\cal A}$, defines a
pure state of ${\cal A}$ [KR; 10.2.5] and so an application of Corollary
9 reveals that $\tilde{\rho}\circ \psi _{I-E}\in {\rm ext}({\cal A}%
_1^{*})$ as required.

If on the other hand $\rho(VE) = 1$, then on arguing as before, $\rho$
defines an extreme point of $({\cal B}_E)^\ast_1$ and annihilates $%
{\cal B}_{I-E}$. We may therefore replace ${\cal B}$ by ${\cal B}_E$
and assume that $I = E$. Since in addition $\ast$-isomorphisms trivially
preserve extreme points of ${\cal B}^\ast_1$, we assume for the moment
that ${\cal B}^{\prime\prime}= {\cal R}_{E_2}$ and that $\psi$ is of
the form

\setcounter{equation}{0} 
\begin{equation}
\psi (A)=E_2W\varphi (A)E_2\quad \mbox{\rm for all}\quad A\in 
{\cal A}.
\end{equation}

We show that it is sufficient to consider only this case. To see this note
that by hypothesis $W$ defines a unitary mapping from $E_1h$ onto $E_2h$,
and hence ${\cal R}_{E_2}\rightarrow {\cal R}_{E_1}:A\rightarrow
E_1W^{*}AWE_1$ defines a spatial $*$-isomorphism from ${\cal R}_{E_2}$
onto ${\cal R}_{E_1}$. Since in addition for all $A\in {\cal R}$ we
have 
\[
E_1\varphi (A)WE_1=W^{*}W\varphi (A)WE_1=E_1W^{*}(E_2W\varphi (A)E_2)WE_1,
\]
it is clear that the one case is $*$-isomorphic to the other, and hence we
may restrict attention to the case outlined in (1) above. Now apply Lemmas 7
and 10 to find a unitary $V\in {\cal R}_{E_2}$ such that 
\[
\omega (A)=\rho (VA)\quad \mbox{\rm for all}\quad A\in {\cal R}_{E_2}
\]
defines a pure state of ${\cal R}_{E_2}$. Denoting the map ${\cal R}%
\rightarrow {\cal R}_{E_2}:A\rightarrow E_2AE_2$ by $\eta _2$, Lemma 1
reveals that $\omega \circ \eta _2$ is a pure state of ${\cal R}_1^{*}$
with 
\[
\rho \circ \eta _2(E_2A)=\rho \circ \eta _2(A)\quad \mbox{\rm for all}
\quad A\in {\cal R}.
\]
Now since $V^{*}V=VV^{*}=E_2$, it can easily be verified that $W^{*}V$ is a
partial isometry with initial projection $E_2$ and final projection $E_1$.
Since then

\[
\begin{array}{lll}
\rho \circ \eta _2(W(W^{*}VA)) & = & \rho \circ \eta _2(E_2VA) \\ 
& = & \rho (VE_2AE_2) \\ 
& = & \omega \circ \eta _2(A)\quad \mbox{\rm for all}\quad
A\in {\cal R},
\end{array}
\]

Corollary 9 reveals that $A \rightarrow \rho \circ \eta_2 (WA)$ defines an
extreme point of ${\cal R}^\ast_1$ which is moreover ultra-weakly
continuous by Lemma 12 and [BR; 2.4.2]. The problem thus reduces to showing
that the adjoint of some Jordan-morphism $\psi : {\cal A} \rightarrow 
{\cal R}$ for which $\psi({\cal A})$ has the stated density condition,
maps ultra-weakly continuous elements of ${\rm ext} ({\cal R}^\ast_1)$
into ${\rm ext}({\cal A}^\ast_1)$. Hence assume $\psi$ to be such a
mapping and let $\rho$ be an ultra-weakly continuous element of ${\rm ext}%
({\cal R}^\ast_1)$. By Lemmas 7 and 10 there exists a unitary $V \in 
{\cal R}$ and an ultra-weakly continuous pure state $\omega$ on ${\cal %
R}$ with $\rho(A) = \omega(V^\ast A)$ and $\rho(VA) = \omega(A)$ for every $%
A \in {\cal R}$.

Now let $(F_\nu )$ be the family of mutually orthogonal central projections
described in the hypothesis. As in the proof of Theorem 5 we may apply [KR;
4.3.14] to conclude that $\omega (F_\nu )\in \{0,1\}$ for every $\nu $, and
then make use of the fact that $\omega (\bigvee F_\nu )=\omega (I)=1$ to
conclude that $\omega (F_\nu )=1$ for precisely one $\nu $, say $\omega
(F_{\nu _0})=1$. Then surely by [KR; 4.3.14],

\begin{equation}
\rho(A) = \omega (V^\ast A) = \omega (V^\ast F_{\nu_0} A) = \rho(F_{\nu_0} A)
\end{equation}

for every $A \in {\cal R}$. Since moreover $\rho(F_{\nu_0} V) = \omega
(F_{\nu_0}) = 1$, it now follows from Lemmas 12 and 13 that the restriction
of $\rho$ to ${\cal R}_{F_{\nu_0}}$ is an ultra-weakly continuous extreme
point of $({\cal R}_{F_{\nu_0}})^\ast_1$. From this fact and (2) it is
clear that we may replace ${\cal R}$ by ${\cal R}_{F_{\nu_0}}$ and
hence effectively assume that $\psi({\cal A})^{\prime\prime}= {\cal R}$
since $(F_{\nu_0} \psi ({\cal A}) F_{\nu_0})^{\prime\prime}= {\cal R}%
_{F_{\nu_0}}$ by hypothesis. Next apply [BR; 3.2.2] to obtain a projection $%
\tilde{E} \in {\cal R} \cap {\cal R}^{\prime}$ such that $\tilde{E}
\psi (A), \ A \in {\cal A}$, defines a $\ast$-homomorphism and $(I - 
\tilde{E}) \psi (A), \ A \in {\cal A}$, a $\ast$-antimorphism. Now by
[KR; 4.3.14] we have $\omega(\tilde{E}) \in \{0,1\}$ and hence again by [KR;
4.3.14], for any $A \in {\cal R}$ we either have

\begin{equation}
\rho (A)=\omega (V^{*}A)\omega (\tilde{E})=\omega (V^{*}\tilde{E}A)=\rho (%
\tilde{E}A)\quad {\rm if}\quad \omega (\tilde{E})=1
\end{equation}

or similarly

\begin{equation}
\rho(A) = \rho((I - \tilde{E}) A) \quad {\rm if} \quad \omega(I - \tilde{E%
} ) = 1.
\end{equation}

We show that we may assume $\psi$ to be either a $\ast$-mor\-phism, or a $\ast$%
-anti\-mor\-phism by con\-sidering the two cases separately.

{\bf Case 1} $(\omega(\tilde{E}) = 1)$ \quad Denote the $C^\ast$-algebra
generated by $\psi(A)$ by ${\cal C}$. Then since ${\cal C}%
^{\prime\prime}= {\cal R}$, we surely have $(\tilde{E}{\cal C} \tilde{E%
} )^{\prime\prime}= {\cal R}$. As $\tilde{E} \psi (\cdot)$ is a $\ast$%
-morphism, $\tilde{E} \psi ({\cal A}) \tilde{E}$ is already a $C^\ast$%
-algebra and hence without too much ado we have $\tilde{E}{\cal C}\tilde{E%
} = \tilde{E} \psi({\cal A})\tilde{E}$. Moreover as $1 = \omega(\tilde{E}%
) = \rho(\tilde{E} V) \leq \Vert \rho |_{\tilde{E}{\cal R}\tilde{E}}
\Vert \leq 1$, it is clear from Lemmas 12 and 13 that the restriction of $%
\rho \circ \eta_{{\cal R}}$ to ${\cal R}_{\tilde{E}}$ is an
ultra-weakly continuous extreme point of $({\cal R}_{\tilde{E}})^\ast_1$.
Since by (2) $\rho$ vanishes on ${\cal R}_{I-\tilde{E}}$, the assertion
follows.

{\bf Case 2} $(\omega(\tilde{E}) = 0$, ie $\omega (I-\tilde{E}) = 1)$
\quad The proof of this case is virtually identical to Case 1, and is
therefore omitted.

Recapitulating, we have thus effectively reduced the situation in question
to 
\[
\psi : {\cal A} \rightarrow \psi({\cal A})^{\prime\prime}= {\cal R}
\]
where ${\cal A}$ is a $C^\ast$-algebra, ${\cal R}$ a von Neumann
algebra, and $\psi$ a $\ast$-(anti)\-mor\-phism with the immediate task at hand
being to show that $\rho \circ \psi \in {\rm ext}({\cal A}^\ast_1)$
for every ultra-weakly continuous extreme point $\rho$ of ${\cal R}%
^\ast_1 $. If now we apply [KR; 10.1.12], then in the notation of [KR;
10.1.12], it follows from for example [BR; 2.4.23] that ${\cal R}$ may be
identified with $P\Phi({\cal C})^{\prime\prime}P$ as far as we are
concerned. Here ${\cal C}$ is of course the $C^\ast$-algebra $\psi(%
{\cal A})$. If next we apply Lemmas 12 and 13, it is clear that $\rho$
extends to an ultra-weakly continuous extreme point of $(\Phi({\cal C}%
)^{\prime\prime})_1$. It follows that we may assume ${\cal C}$ (ie $\psi(%
{\cal A}))$ to be in its universal representation. Moreover a slight
modification of [St\o 3, Lemma 3.1] now reveals that we may extend $\psi$ to
an ultra-weakly continuous $\ast$-(anti)morphism $\psi^{\ast\ast} : {\cal %
A}^{\ast\ast} \rightarrow {\cal R}$. Since $\rho$ is ultra-weakly
continuous, it now follows from the Kaplansky density theorem [BR; 2.4.16]
and the fact that $\rho$ is trivially continuous under a slightly stronger
topology on ${\cal R}$ (the $\sigma$-strong$^{\ast}$ topology [BR;
p270]), that $\rho$ assumes its norm on $\psi({\cal A})_1$ and hence also
on the larger set $\psi^{\ast\ast}({\cal A}^{\ast\ast})_1$. However as $%
\psi^{\ast\ast}$ is a $\ast$-(anti)morphism, we have that $\Vert
\psi^{\ast\ast} \Vert = 1$ and $\psi^{\ast\ast}(({\cal A}^{\ast\ast})_1)
= \psi^{\ast\ast}({\cal A}^{\ast\ast})_1$. Hence we may select a sequence 
$(A_n) \subset ({\cal A}^{\ast\ast})_1$ with $\rho \circ \psi^{\ast\ast}
(A_n) \rightarrow \Vert \rho \Vert = 1$. Since $\Vert \rho \circ
\psi^{\ast\ast} \Vert \leq \Vert \rho \Vert \Vert \psi^{\ast\ast} \Vert = 1$%
, it is clear that $\rho \circ \psi^{\ast\ast}$ is then a norm-one
ultra-weakly continuous functional on ${\cal A}$. Thus by [KR; 7.3.2]
there exists a normal state $\omega_0$ on ${\cal A}^{\ast\ast}$ and a
partial isometry $W \in {\cal A}^{\ast\ast}$ with

\begin{equation}
\omega_0 (A) = \rho \circ \psi^{\ast\ast} (WA) \quad {\rm and} \quad
\omega_0 (W^\ast A) = \rho \circ \psi^{\ast\ast} (A)
\end{equation}

for all $A \in {\cal A}^{\ast\ast}$. Since $\psi^{\ast\ast}$ is a $\ast$%
-(anti)morphism, $\psi^{\ast\ast}(W)$ is a partial isometry in ${\cal R}
= \psi({\cal A})^{\ast\ast}$. We conclude by considering two cases.

{\bf Case 1} \ $(\psi^{\ast\ast}$ a $\ast$-homomorphism) \quad If we
define $\omega_1 : \psi({\cal A}) \rightarrow {\bf C}$ by 
\[
\omega_1 (\psi(A)) = \rho(\psi^{\ast\ast}(W) \psi(A)) = \rho
(\psi^{\ast\ast} (WA)) = \omega_0 (A)
\]
for every $A \in {\cal A}$ and uniquely extend $\omega_1$ to $\psi(%
{\cal A})^{\ast\ast} = {\cal R},$ then since $\omega_1(I) = \omega_1
(\psi^{\ast\ast}(I)) = \omega_0(I) = 1,$ \ $\omega_1$ is a state of $%
{\cal R}$ [KR; 4.3.2], which is moreover ultra-weakly continuous [KR;
10.1.1]. Lemma 10 and the fact that $\omega_1$ is ``related'' to $\rho$ by
means of $\psi^{\ast\ast}(W)$ by (5) above, now reveals that $\omega_1$ is
in fact pure. But then $\omega_1 \circ \psi$ must be a pure state of $%
{\cal A}$ by Theorem 5, and hence since $\omega_1 \circ \psi$ is just the
restriction of $\omega_0$ to ${\cal A}$, \ $\omega_0$ is a (normal) pure
state by Lemma 4. Thus by Corollary 9 and (5), the restriction of $\rho
\circ \psi^{\ast\ast}$ to ${\cal A}$, \ $\rho \circ \psi$, is an extreme
point of ${\cal A}^\ast_1$.

{\bf Case 2} $(\psi ^{**}$ an anti-morphism) \quad As a start we observe
that since $\rho $ is an ultra-weakly continuous extreme point of ${\cal R%
}_1^{*}$, it is an easy exercise to show that the same is true of $\rho ^{*}$%
. Now define $\omega _1:\psi ({\cal A})\rightarrow {\bf C}$ by 
\[
\omega _1(\psi (A))=\rho ^{*}(\psi ^{**}(W^{*})\psi (A))\quad 
    \mbox{\rm for all}\quad A\in {\cal A}.
\]
Observe that for any self-adjoint $A\in {\cal A}$ we have

\begin{eqnarray}
\omega_1 (\psi(A)) & = & \rho^\ast (\psi^{\ast\ast} (W^\ast) \psi(A)) = 
\overline{\rho((\psi^{\ast\ast}(W^\ast) \psi(A))^\ast)}  \nonumber\\
& = & \overline{\rho(\psi(A) \psi^{\ast\ast}(W))} = \overline{%
\rho(\psi^{\ast\ast}(WA))} \\
& = & \overline{\omega_0(A)} = \omega_0(A).  \nonumber
\end{eqnarray}

By linearity the above holds for ${\rm span}({\cal A}_{sa}) = {\cal %
A}$. Since $\omega_1(I) = \omega_1(\psi(I)) = \omega_0(I) = 1$ with $\Vert
\omega_1 \Vert \leq \Vert \rho^\ast \Vert \Vert \psi^{\ast\ast} \Vert = 1, \
\omega_1 $ is a state of $\psi({\cal A})$. By Lemma 10 applied to the
fact that $\omega_1$ is ``related'' to $\rho^\ast$ by means of $%
\psi^{\ast\ast}(W^\ast),$ this state is pure, and by Theorem 5 $\omega_1
\circ \psi$ is then a pure state of ${\cal A}$. But from (6) above, $%
\omega_1 \circ \psi = \omega_0 |_{{\cal A}}$. Thus by Lemma 4, $\omega_0$
is the unique ultra-weakly continuous pure extension of $\omega_1 \circ \psi$
to all of ${\cal A}^{\ast\ast}$. As in case 1, Corollary 9 now reveals
that $\rho \circ \psi$ is an extreme point of ${\cal A}^\ast_1.$

It remains to verify that (c) follows from (b).

{\bf (b) $\Longrightarrow $ (c)} \quad To see this assume (b) to be true,
and apply Zorn's lemma to obtain a maximal set ${\cal M}$ of pure states
of ${\cal B}$ for which the associated (irreducible) GNS representations
of ${\cal B}$ are pairwise inequivalent. Then surely $\pi
=\bigoplus_{\omega \in {\cal M}}\pi _\omega $ corresponds to the reduced
atomic representation of ${\cal B}$ (here $(\pi _\omega ,h_\omega ,\Omega
_\omega )$ is the canonical representation engendered by $\omega \in 
{\cal M}$). Now partition ${\cal M}$ into the two disjoint classes $%
{\cal M}_1$ and ${\cal M}_2$ where for $i=1,2,\ \omega \in {\cal M}%
_1$ (respectively $\omega \in {\cal M}_2)$ if and only if $\omega \in 
{\cal M}$ and $\pi _\omega \circ \psi $ is of the form (1) (respectively
(2)) as described in (b). Now let $E\in B(\bigoplus_{\omega \in {\cal M}%
}h_\omega )$ be the canonical projection of $\bigoplus_{\omega \in {\cal M%
}}h_\omega $ onto the subspace corresponding to 
\[
\bigoplus_{\omega \in {\cal M}_1}h_\omega .
\]

By construction it is now clear that $E\in (\bigoplus_{\omega \in {\cal M}%
}\pi _\omega ({\cal B}))^{\prime }$. Moreover since on restriction to the
subspace corresponding to $h_\omega $ (for some given $\omega \in {\cal M}%
)$ \ $E$ corresponds to either the identity or the zero-operator on $%
h_\omega $, it is clear from [KR; 10.3.10] that $E\in \pi ({\cal B})$
since ${\cal B}$ is unital. Now given any irreducible representation $%
\alpha $ of $\pi ({\cal B})_{I-E}$, it follows from the maximality of $%
{\cal M}$ that $\alpha (\pi ({\cal B})_{I-E})$ is equivalent to $\pi
_{\omega _0}({\cal B})$ for some $\omega _0\in {\cal M}$. Hence $%
A\rightarrow \alpha ((I-E)\pi (A)),\ A\in {\cal B}$, has the same kernel
as $\pi _{\omega _0}$. However the presence of $I-E$ leads us to conclude
that for any $\omega \in {\cal M}_1$ there is some $A\in {\cal B}$
with $\pi _\omega (A)\neq 0$ and $\alpha ((I-E)\pi (A))=0$. Hence $\omega
_0\in {\cal M}_2$. From this fact and the equivalence of $\alpha (\pi (%
{\cal B})_{I-E})$ and $\pi _{\omega _0}({\cal B})$, we now conclude
that $\alpha \circ ((I-E)\pi \circ \psi (\cdot ))$ is of the form described
in (b(2)) as required. To conclude the proof we consider the mapping $%
{\cal A}\rightarrow \pi ({\cal B})_E:A\rightarrow E(\pi \circ \psi
(A))E$ and show that it is of the requisite form described in (c(2)). Recall
that each $\pi _\omega \circ \psi ,\ \omega \in {\cal M}_1$, is of the
form

\begin{equation}
\pi _\omega \circ \psi (A)=V_\omega ^{*}\alpha _\omega (A)U_\omega \quad 
\mbox{\rm for all}\quad A\in {\cal A}
\end{equation}

for some irreducible $*$-(anti)morphism $\alpha _\omega $ from ${\cal A}$
into $B(k_\omega )$. Now let $h=\bigoplus_{{\cal M}_1}k_\omega $ and $%
{\cal R}=\bigoplus_{{\cal M}_1}B(k_\omega )$. If for each $\omega \in 
{\cal M}_1$ we let $F_\omega $ be the canonical projection of $h$ onto
the subspace corresponding to $k_\omega $, and if we let $\varphi
=\bigoplus_{{\cal M}_1}\alpha _\omega $, then it may easily be verified
that $\varphi $ is a Jordan $*$-morphism (since each $\alpha _\omega $ is),
that $(F_\omega )_{{\cal M}_1}\subset {\cal R}\cap {\cal R}^{\prime
\prime }$ is mutually orthogonal with $\sum_{{\cal M}_1}F_\omega =I$, and
that $\varphi ({\cal A})=\bigoplus_{{\cal M}_1}\alpha _\omega (%
{\cal A})$ has the required density property in terms of $(F_\omega )_{%
{\cal M}_1}$. Finally note that since by hypothesis the mappings $%
V_\omega :h_\omega \rightarrow k_\omega $ and $U_\omega :h_\omega
\rightarrow k_\omega $, \ $\omega \in {\cal M}_1$, referred to in (7) are
injective partial isometries, it follows that the same is true of 
\[
V:\bigoplus_{{\cal M}_1}h_\omega \rightarrow h,\ U:\bigoplus_{{\cal M}%
_1}h_\omega \rightarrow h\quad {\rm where}\quad V=\bigoplus_{{\cal M}%
_1}V_\omega
\]
and $U=\bigoplus_{{\cal M}_1}U_\omega $. Since effectively $\bigoplus_{%
{\cal M}_1}h_\omega $ appears as the image of $E$, we may suppose $V$ and 
$U$ to be acting from $E(\bigoplus_{{\cal M}}h_\omega )$. Thus by
construction 
\[
\pi \circ \psi _E(A)=V^{*}\varphi (A)U\quad \mbox{\rm for all}
\quad A\in {\cal A}.
\]
The injectivity of $U$ and $V$ now imply that $U^{*}U=I_h=V^{*}V$, with in
addition $W=UV^{*}=\bigoplus_{{\cal M}_1}U_\omega V_\omega ^{*}\in
\bigoplus_{{\cal M}_1}B(k_\omega )={\cal R}$. Notice further that now 
\[
W^{*}W=V(U^{*}U)V^{*}=R_V\quad {\rm and}\quad WW^{*}=U(V^{*}V)U^{*}=R_U.
\]
Hence $W$ is a partial isometry with initial projection $R_V$, the range
projection of $V$, and final projection $R_U$, the range projection of $U$.
Let $E_1=R_V$ and $E_2=R_U$. Then 
\[
\pi \circ \psi _E(A)=V^{*}\varphi (A)U=U^{*}W\varphi (A)U
\]
for all $A\in {\cal A}$. However note that since $U$ is an injective
partial isometry from $E(\bigoplus_{{\cal M}}h_\omega )$ onto $E_2(h)$,
the mapping $\pi ({\cal B})_E\rightarrow {\cal R}_{E_2}:A\rightarrow
UAU^{*}$ now turns out to be a spatial $*$-isomorphism. Hence up to a $*$%
-isomorphism $(\pi ({\cal B})_E)^{\prime \prime }$ appears as ${\cal R}%
_{E_2}$ with $\pi \circ \psi _E$ corresponding to the map 
\[
U(U^{*}W\varphi (A)U)U^{*}=E_2W\varphi (A)E_2,\qquad A\in {\cal A}.
\]
In a similar fashion one can show that up to a spatial $*$-isomorphism
induced by $V,$ $(\pi ({\cal B})_E)^{\prime \prime }$ appears as $%
{\cal R}_{E_1}$, with $\pi \circ \psi _E$ now corresponding to the map 
\[
E_1\varphi (A)WE_1,\qquad A\in {\cal A}.
\]
\hfill $\Box $

With Theorems 5 and 19 now at our disposal, the conclusion regarding maps
with pure state preserving adjoints is now immediately obvious.

{\bf Corollary 20} \quad Let ${\cal A}, {\cal B}$ be $C^\ast$%
-algebras and $\psi : {\cal A} \rightarrow {\cal B}$ a linear map with
the property that $\omega \circ \psi \in {\cal P}_{{\cal A}}$ whenever 
$\omega \in {\cal P}_{{\cal B}}$. Then $\rho \circ \psi \in {\rm ext%
} ({\cal A}^\ast_1)$ whenever $\rho \in {\rm ext}({\cal B}^\ast_1)$.

{\bf Proof} \quad Consider Theorem 5 alongside Theorem 19. \hfill $\Box$

\section{Maps on some specific spaces}

In conclusion, to provide information about what may happen in the
commutative case, we consider the action of maps with ``extreme point
preserving'' adjoints on some uniform algebras, hereby considering the
commutative $C^{*}$-algebras in particular, where the maps in question
reduce to ordinary compositions of a multiplication and of a composition
operator.

If $A\subset C(X)$ is a uniform algebra over $X$, let $\partial _A$ be its $%
\check{S}$ilov boundary. Denote by $p_A$ the so-called Choquet (or strong)
boundary of $A$ (this is the set of all $p$-points in $X$, or generalized
peak points --- see [Gam]). Also, if $x\in X$ and $\delta _x$ is the
functional corresponding to the evaluation at $x$, denote by $\eta _x$ its
restriction to $A$. The following Lemma is very probably well-known, but we
have not been able to find a reference for it:

{\bf Lemma 21} \quad Let $A\subset C(X)$ be a uniform algebra. Then 
\[
{\rm  ext}(A_1^{*})=\partial D\cdot p_A\qquad {\rm and}\qquad \overline{%
{\rm ext}(A_1^{*})}=\partial D\cdot \partial _A\quad .
\]

{\bf Proof} \quad ``$\supset $'' : Let $x\in X$ be a $p$-point, and
suppose we have $0\leq t\leq 1$ and $\varphi ,\psi \in B_{A^{*}}$ such that $%
\eta _x=t\varphi +(1-t)\psi $. Clearly, if $\tilde{\varphi }$ (resp. $\tilde{%
\psi}$) are Hahn-Banach extensions of $\varphi $ (resp. $\psi $), then $\xi
:=t\tilde{\varphi }+(1-t)\tilde{\psi}$ is an extension of $\eta _x$. But
since for $p$-points such an extension is unique [G], we must have $\xi
=\delta _x$. Now, $\delta _x\in $ext $(C(X)_1^{*})$, and so $\delta _x=%
\tilde{\varphi }=\tilde{\psi}$ and, in particular, $\eta _x=\varphi =\psi $.
Consequently, $\eta _x\in $ext $(A_1^{*})$.

``$\subset $'' : We first prove that ext $(A_1^{*})\subset \partial D\cdot
\partial _A$. Let $\delta \in $ext $(A_1^{*})$, and consider the set $%
S:=\{\mu \in C(X)^{*}\,:\left\| \,\mu \right\| =1,\int f\,d\mu
=\delta (f),\forall \,f\in A\,\,{\rm and}\,\,{\rm supp}\,\mu \subset
\partial _A\}$. $S$ is w$^{*}$-compact in $C(X)^{*}$ and can be
partially ordered by setting $\mu \preceq \nu $ iff supp $\mu \subset $ supp 
$\nu $. By the w$^{*}$-compactness of $S$ Zorn's lemma applies to
give a $\mu \in S$ with minimal support $K\subset \partial _A$. If 
$K$ is a point $x$, then clearly $\delta (f)=e^{it}\delta _x(f)$ for some
real $t$, and we are done. If we allow two distinct $x,y$ in $K$ we can
derive a contradiction as follows. Choose an open neighborhood $E$ of $x$
such that $y\notin \overline{E}$. Then $\left| \mu \right| (E)\neq 0\neq
\left| \mu \right| (X\setminus E)$ and we may write 
\[
\delta (f)=\left| \mu \right| (E)\left( \frac 1{\left| \mu \right|
(E)}\int_Ef\,d\mu \right) +(1-\left| \mu \right| (E))\left( \frac 1{1-\left|
\mu \right| (E)}\int_{X\setminus E}f\,d\mu \right) 
\]
for all $f\in A$. Since $\delta $ is in ext $(A_1^{*})$, we must have 
\[
\delta (f)=\frac 1{\left| \mu \right| (E)}\int_Ef\,d\mu\qquad 
\mbox{\rm for all}\quad\,f\in A.
\]
Hence, $1_E\mu /\left| \mu \right| (E)$ is in $S$ but its support $%
\overline{E}$ is strictly smaller than the minimal $K$.

Suppose now $x\in \partial _A\setminus p_A$. Then, by the Bishop-deLeeuw
theorem [Gam, 12.9] there exists a probability measure $\mu $ on the $\sigma 
$-algebra generated by $p_A$ and the Borel sets in $X$, such that $\mu
(p_A)=1$ and $\delta _x(f)=f(x)=\int f\,d\mu $ for all $f\in A$. Let $y\in
p_A\cap {\rm supp}\mu $, and let $f\in A$ be a function with $f(y)=1$ and 
$\left| f(x)\right| <1/2$. Choosing a small open neighborhood $E$ of $y$ we
can ensure that $\left| {\frac 1{\mu (E)}\int_Ef\,d\mu \,-\,1}\right| <1/2$
and $\mu (p_A\setminus E)\neq 0$ (clearly, $\mu $ cannot be concentrated
on $y$). So, we write 
\[
\delta _x(g)={\mu }(E)\left( \frac 1{{\mu }(E)}\int_Eg\,d\mu \right) +(1-{%
\mu }(E))\left( \frac 1{1-{\mu }(E)}\int_{E^c}g\,d\mu \right) 
\]
for all $g\in A$. If $\delta _x$ were an extreme point of $B_{A^{*}}$, then
we would have for the above $f$ the contradiction $\delta _x(f)=\frac 1{\mu
(E)}\int_Ef\,d\mu $.

Finally, the statement about $\overline{{\rm ext}(A_1^{*})}$ follows
immediately from the fact $\overline{p_A}=\partial _A$.\hfill$\Box $

{\bf Corollary 22} \quad If $A$ is logmodular then ext $%
(A_1^{*})=\partial D\cdot \partial _A.$

This follows from the fact that multiplicative functionals on logmodular
algebras have unique representing measures [Gam,II.4.2], and from
[Gam,II.11.3] which states that a point $x\in\partial_A$ is a $p$-point 
if and only if
the point mass at $x$ is the only representing measure for the point
evaluation at $x$.

Using the above corollary (and classical results such as: every point in 
$partial D$ is a peak point for the disc algebra, $H_\infty $ is logmodular
--- see [Gam, Gar]) it is now more or less immediate to deduce the following
(where for the sake of simplicity we call extremal
an operator $\Phi $ such that $\Phi ^{*}$
sends extreme points of the domain ball to extreme points of the range ball):

{\bf Theorem 23} \quad (a) An operator $\Phi $ on the disc algebra 
${\cal A}$ is extremal if and only if $\Phi =M_\psi C_\varphi $, where $\psi$ 
and $\varphi$ are finite Blaschke products. (b) An operator $\Phi$ on 
$H_\infty$ is extremal if and only if $\Phi =M_\psi C_\varphi$, where $\psi$
and $\varphi$ are inner functions.\hfill$\Box$

\vspace{3cm}

The authors would like to thank prof H. Jarchow for initiating the contact
between them, the result of which is this paper.

\end{document}